\documentclass[a4]{article}

\font\es=eufm10

\def\gge{\mbox{\es {e}}} 
\def\gf{\mbox{\es {f}}} 
\def\gg{\mbox{\es {g}}}

\def\gC{\mbox{\es {C}}}

\def\gJ{\mbox{\es {J}}}

\def\gP{\mbox{\es {P}}}

\def\so{\mbox{\es {so}}}

\def\spin{\mbox{\es {spin}}}
\def\tr{\mbox{\rm {tr}}}

\def\diag{\mbox{\rm {diag}}}
\def\Iso{\mbox{\rm {Iso}}}

\def\Ker{\mbox{\rm {Ker}}}
\def\ad{\mbox{\rm {ad}}}
\def\ov{\overline}
\def\wti{\widetilde}

\def\dfrac#1#2{\displaystyle \frac{#1}{#2}}

\def\H{\mbox{\boldmath $H$}}

\def\R{\mbox{\boldmath $R$}}

\def\Z{\mbox{\boldmath $Z$}} 
 
\def\sR{\mbox{\boldmath $\scriptstyle{R}$}} 
 
\def\0{\mbox{\boldmath {0}}}    
\def\1{\mbox{\boldmath {1}}}      
\def\2{\mbox{\boldmath {2}}}      
\def\3{\mbox{\boldmath {3}}}      
\def\4{\mbox{\boldmath {4}}}      
\def\5{\mbox{\boldmath {5}}}      
\def\6{\mbox{\boldmath {6}}}      
\def\7{\mbox{\boldmath {7}}}      
\def\8{\mbox{\boldmath {8}}}      
\def\9{\mbox{\boldmath {9}}}

\begin{document}
\baselineskip=14pt

\begin{center}
{\Large{\bf Decomposition of spinor groups by the involution $\sigma'$}}
\end{center}
\begin{center}
{\Large{\bf  in exceptional Lie groups}}
\end{center}
\vspace{5mm}
\begin{center}
{\large{Toshikazu Miyashita}}
\end{center}
\vspace{6mm}

{\large{\bf Introduction}}
\vspace{3mm}

The compact exceptional Lie groups $F_4, E_6, E_7$ and $E_8$ have spinor groups as a subgroup as follows.
$$
\begin{array}{l}
   F_4 \supset Spin(9) \supset Spin(8) \supset Spin(7) \supset \cdots \supset Spin(1) \ni 1 
\vspace{-0.5mm}\\
   \;\cap 
\vspace{-0.5mm}\\
   E_6 \supset Spin(10)
\vspace{-0.5mm}\\
\;\cap
\vspace{-0.5mm}\\
   E_7 \supset Spin(12) \supset Spin(11)
\vspace{-0.5mm}\\
\;\cap 
\vspace{-0.5mm}\\
   E_8 \supset Ss(16) \supset Spin(15) \supset Spin(14) \supset Spin(13)
\end{array}$$
On the other hand, we know the involution $\sigma'$ induced an element $\sigma' \in Spin(8) \subset F_4 \subset E_6 \subset E_7 \subset E_8$. Now, in this paper, we determine the group structures of $(Spin(n))^{\sigma'}$ which are the fixed subgroups by the involution $\sigma'$. Our results are as follows.
$$
\begin{array}{ll}
   F_4 & \quad (Spin(9))^{\sigma'} \cong Spin(8)
\vspace{1mm}\\  
   E_6 & \quad (Spin(10))^{\sigma'} \cong (Spin(2) \times Spin(8))/\Z_2
\vspace{1mm}\\  
   E_7 & \quad (Spin(11))^{\sigma'} \cong (Spin(3) \times Spin(8))/\Z_2   
\vspace{1mm}\\  
   {} & \quad (Spin(12))^{\sigma'} \cong (Spin(4) \times Spin(8))/\Z_2   
\vspace{1mm}\\  
   E_8 & \quad (Spin(13))^{\sigma'} \cong (Spin(5) \times Spin(8))/\Z_2
\vspace{1mm}\\  
   {} & \quad (Spin(14))^{\sigma'} \cong (Spin(6) \times Spin(8))/\Z_2
\end{array}$$
Needless to say, the spinor groups appeared in the first term have relation
$$
   Spin(2) \subset Spin(3) \subset Spin(4) \subset Spin(5) \subset Spin(6). $$
One of our aims is to find these groups explicitly in the exceptional groups. In the group $E_8$, we conjecture that
$$
   (Spin(15))^{\sigma'} \cong (Spin(7) \times Spin(8))/\Z_2, \;
   (Ss(16))^{\sigma'} \cong (Spin(8) \times Spin(8))/(\Z_2 \times \Z_2), $$
however, we can not realize explicitly.
\vspace{1mm}

This paper is closely in connection with the preceding papers [2],[3],[4] and may be a continuation of [2],[3],[4] in some sense. 
\vspace{4mm}

\begin{center}
{\large{\bf 1. Group $F_4$}}
\end{center}
\vspace{3mm}

We use the same notation as in [5] (however, some will be rewritten). For example,
\vspace{-4mm}

the Cayley algebra $\gC = \H \oplus \H e_4$,
\vspace{1mm}

the exceptional Jordan algebra $\gJ = \{X \in M(3, \gC) \, | \, X^* = X \}$, the Jordan multiplication $X \circ Y$, the inner product $(X, Y)$ and the elements $E_1, E_2, E_3 \in \gJ$,
\vspace{1mm}

the group $F_4 = \{\alpha \in \Iso_{\sR}(\gJ) \, | \, \alpha(X \circ Y) = \alpha X \circ \alpha Y \}$, and the element $\sigma \in F_4: \sigma X = DXD, D = \diag(1, -1, -1), X \in \gJ$ and the element $\sigma' \in F_4: \sigma'X = D'XD', D' = \diag(-1, -1, 1), X \in \gJ$,
\vspace{1mm}

the groups $SO(8) = SO(\gC)$ and $Spin(8) = \{(\alpha_1, \alpha_2, \alpha_3) \in SO(8) \times SO(8) \times SO(8) \, | \, (\alpha_1x)(\alpha_2y) = \ov{\alpha_3\ov{(xy)}} \}$.
\vspace{3mm}

{\bf Proposition 1.1.} \qquad \qquad $(F_4)_{E_1} \cong Spin(9)$.
\vspace{2mm}

{\bf Proof.} We define a 9 dimensional $\R$-vector space $V^9$ by
$$
   V^9 = \{X \in \gJ \, | \, E_1 \circ X = 0, \tr(X) = 0 \}
          = \Big\{ \pmatrix{0 & 0 & 0 \cr
                            0 & \xi & x \cr
                            0 & \ov{x} & -\xi} \Big| \, \xi \in \R, x \in \gC \Big\}$$
with the norm $1/2(X, X) = \xi^2 + \ov{x}x$. Let $SO(9) = SO(V^9)$. Then, we have $(F_4)_{E_1}/\Z_2 
\vspace{0.5mm}
\cong SO(9)$, $\Z_2 = \{1, \sigma \}$. Therefore, $(F_4)_{E_1}$ is isomorphic to $Spin(9)$ as a double covering group of $SO(9)$. (In details, see [5],[8].)
\vspace{2mm}

Now, we shall determine the group structure of $(Spin(9))^{\sigma'}$.
\vspace{3mm}

{\bf Theorem 1.2.} \qquad \qquad $(Spin(9))^{\sigma'} \cong Spin(8)$.
\vspace{2mm}

{\bf Proof.} Let $Spin(9) = (F_4)_{E_1}$. Then, the map $\varphi_1 : Spin(8) \to (Spin(9))^{\sigma'}$,
$$
   \varphi_1(\alpha_1, \alpha_2, \alpha_3)X = 
     \pmatrix{\xi_1 & \alpha_3x_3 & \ov{\alpha_2x_2} \cr
              \ov{\alpha_3x_3} & \xi_2 & \alpha_1x_1 \cr
              \alpha_2x_2 & \ov{\alpha_1x_1} & \xi_3}, \quad X \in \gJ $$
gives an isomorphism as groups. (In details, see [3].)
\vspace{4mm}

\begin{center}
{\large{\bf 2. Group $E_6$}}
\end{center}
\vspace{3mm}

We use the same notation as in [5] (however, some will be rewritten). For example,
\vspace{-4mm}

the complex exceptional Jordan algebra $\gJ^C = \{X \in M(3, \gC^C) \, | \, X^* = X \}$, the Freudenthal multiplication $X \times Y$ and the Hermitian inner product $\langle X, Y \rangle$,
\vspace{1mm}

the group $E_6 = \{\alpha \in \Iso_C(\gJ^C) \, | \, \alpha X \times \alpha Y = \tau\alpha\tau(X \times Y), \langle \alpha X, \alpha Y \rangle = \langle X, Y \rangle \}$, and the natural inclusion $F_4 \subset E_6$,
\vspace{1mm}

any element $\phi$ of the Lie algebra $\gge_6$ of the group $E_6$ is uniquely expressed as $\phi = \delta + i\wti{T}, \delta \in \gf_4, T \in \gJ_0$, where $\gJ_0 = \{T \in \gJ \, | \, \tr(T) = 0 \}$.
\vspace{3mm}

{\bf Proposition 2.1.} \qquad \qquad $(E_6)_{E_1} \cong Spin(10)$.
\vspace{2mm}

{\bf Proof.} We define a 10 dimensional $\R$-vector space $V^{10}$ by
$$
   V^{10} = \{X \in \gJ^C \, | \, 2E_1 \times X = - \tau X \}
          = \Big\{ \pmatrix{0 & 0 & 0 \cr
                            0 & \xi & x \cr
                            0 & \ov{x} & -\tau\xi} \Big| \, \xi \in C, x \in \gC \Big\}$$
with the norm $1/2\langle X, X \rangle = (\tau\xi)\xi + \ov{x}x$. Let $SO(10) = SO(V^{10})$. Then, we 
\vspace{0.5mm}
have $(E_6)_{E_1}/ \Z_2 \cong SO(10)$, $\Z_2 = \{1, \sigma \}$. Therefore, $(E_6)_{E_1}$ is isomorphic to $Spin(10)$ as a double covering group of $SO(10)$. (In details, see [5],[8].)
\vspace{3mm}

{\bf Lemma 2.2.} {\it For $\nu \in Spin(2) = U(1) = \{\nu \in C \, | \, (\tau\nu)\nu = 1 \}$, we define a $C$-linear transformation $\phi_1(\nu)$ of $\gJ^C$ by
$$
   \phi_1(\nu)X 
   = \pmatrix{\xi_1 & \nu x_3 & \nu^{-1}\ov{x}_2 \cr
              \nu\ov{x}_3 & \nu^2\xi_2 & x_1 \cr
              \nu^{-1}x_2 & \ov{x}_1 & \nu^{-2}\xi_3}, \quad X \in \gJ^C. $$
Then, $\phi_1(\nu) \in ((E_6)_{E_1})^{\sigma'}$.}
\vspace{3mm}

{\bf Lemma 2.3.} {\it Any element $\phi$ of the Lie algebra $((\gge_6)_{E_1})^{\sigma'}$ of the group $((E_6)_{E_1})^{\sigma'}$ is expressed by
$$
    \phi = \delta + it(E_2 - E_3)^\sim, \quad \delta \in ((\gf_4)_{E_1})^{\sigma'} = \so(8),\,\, t \in \R. $$
In particular, we have}
$$
        \dim(((\gge_6)_{E_1})^{\sigma'}) = 28 + 1 = 29. 
\vspace{2mm}$$

Now, we shall determine the group structure of $(Spin(10))^{\sigma'}$.
\vspace{3mm}

{\bf Theorem 2.4.} \quad $(Spin(10))^{\sigma'} \cong (Spin(2) \times Spin(8))/\Z_2, \Z_2 = \{(1, 1), (-1, \sigma) \}$.
\vspace{2mm}

{\bf Proof.} Let $Spin(10) = (E_6)_{E_1}$, $Spin(2) = U(1) \subset ((E_6)_{E_1})^{\sigma'}$(Lemma 2.2) and $Spin(8) = ((F_4)_{E_1})^{\sigma'} \subset ((E_6)_{E_1})^{\sigma'}$ \,\,(Theorem 1.2, Proposition 2.1). Now, we define a map $\varphi : Spin(2) \times Spin(8) \to (Spin(10))^{\sigma'}$ by
$$
       \varphi(\nu, \beta) = \phi_1(\nu)\beta. $$
Then, $\varphi$ is well-defined : $\varphi(\nu, \beta) \in (Spin(10))^{\sigma'}$. Since $\phi_1(\nu)$ and $\beta$ are commutative, $\varphi$ is a homomorphism. $\Ker\,\varphi = \{(1, 1), (-1, \sigma) \}$. Since $(Spin(10))^{\sigma'}$ is connected and $\dim(\spin(2) \oplus \spin(8)) = 1 + 28 = 29 = \dim((\spin(10)^{\sigma'}))$ (Lemma 2.3), $\varphi$ is onto. Thus, we have the isomorphism $(Spin(2) \times Spin(8))/\Z_2 \cong (Spin(10))^{\sigma'}$.
\vspace{2mm}

\begin{center}
{\large{\bf 3. Group $E_7$}}
\end{center}
\vspace{3mm}

We use the same notation as in [6](however, some will be rewritten). For example,
\vspace{-4mm}

the Freudenthal $C$-vector space $\gP^C = \gJ^C \oplus \gJ^C \oplus C \oplus C$, the Hermitian inner product $\langle P, Q \rangle$,
\vspace{1mm}

for $P, Q \in \gP^C$, the $C$-linear map $P \times Q$ : $\gP^C \to \gP^C$,
\vspace{1mm}

the group $E_7 = \{\alpha \in \Iso_C(\gP^C) \, | \, \alpha(X \times Y)\alpha^{-1} = \alpha P \times \alpha Q, \langle \alpha P, \alpha Q \rangle = \langle P, Q \rangle \}$, the natural inclusion $E_6 \subset E_7$ and elements $\sigma, \sigma' \in F_4 \subset E_6 \subset E_7$, $\lambda \in E_7$,
\vspace{1mm}

any element ${\mit\Phi}$ of the Lie algebra $\gge_7$ of the group $E_7$ is uniquely expressed as ${\mit\Phi} = {\mit\Phi}(\phi, A, -\tau A, \nu), \phi \in \gge_6, A \in \gJ^C, \nu \in i\R$.
\vspace{2mm}

In the following, the group $((Spin(10))^{\sigma'})_{F_1(x)}$ is defined by
$$
   ((Spin(10))^{\sigma'})_{F_1(x)} = \{\alpha \in (Spin(10))^{\sigma'} \, | \, \alpha F_1(x) = F_1(x) \;\;\mbox{{\it for all}}\; x \in \gC \}, $$
where                        $F_1(x) = \pmatrix{0 & 0 & 0 \cr
                                                0 & 0 & x \cr
                                                0 & \ov{x} & 0} \in \gJ$.
\vspace{3mm}

{\bf Proposition 3.1.} \qquad \quad $((Spin(10))^{\sigma'})_{F_1(x)} \cong Spin(2)$.
\vspace{2mm}                       

{\bf Proof.} Let $Spin(10) = (E_6)_{E_1}$ and $Spin(2) = U(1) = \{\nu \in C \, | \, (\tau\nu)\nu = 1 \}$. We consider the map $\phi_1 : Spin(2) \to ((Spin(10))^{\sigma'})_{F_1(x)}$ defined in Section 2. Then, $\phi_1$ is well-defined : $\phi_1(\nu) \in ((Spin(10))^{\sigma'})_{F_1(x)}$. We shall show that $\phi_1$ is onto. From $((Spin(10))^{\sigma'})_{F_1(x)} \subset (Spin(10))^{\sigma'}$, we see that for $\alpha \in ((Spin(10))^{\sigma'})_{F_1(x)}$, there exist $\nu \in Spin(2)$ and $\beta \in Spin(8)$ such that $\alpha = \varphi(\nu, \beta)$ (Theorem 2.4). Further, from $\alpha F_1(x) = F_1(x)$ and $\phi_1(\nu)F_1(x) = F_1(x)$, we have $\beta F_1(x) = F_1(x)$. Hence, $\beta = (1, 1, 1)$ or $(1, -1, -1) = \sigma$  by the principle of triality. Hence, $\alpha =\phi_1(\nu)$ or $\phi_1(\nu)\sigma$. However, in the latter case, from $\sigma = \phi_1(-1)$, we have $\alpha = \phi_1(\nu)\phi_1(-1) = \phi_1(-\nu)$. Therefore, $\phi_1$ is onto. $\Ker\,\phi_1 = \{1\}$. Thus, we have the isomorphism $Spin(2) \cong ((Spin(10))^{\sigma'})_{F_1(x)}$.
\vspace{2mm}

We define $C$-linear maps $\kappa, \mu : \gP^C \to \gP^C$ respectively by
$$
\begin{array}{l}
   \kappa(X, Y, \xi, \eta) = (-\kappa_1X, \kappa_1Y, -\xi, \eta), \;\; \kappa_1 X = (E_1, X)E_1 - 4E_1 \times (E_1 \times X),
\vspace{1mm}\\
   \mu(X, Y, \xi, \eta) = (2E_1 \times Y + \eta E_1, 2E_1 \times X + \xi E_1, (E_1, Y), (E_1, X)).
\end{array}$$
Their explicit forms are
$$
\begin{array}{l}
   \kappa(X, Y, \xi, \eta) = \Big(\pmatrix{-\xi_1 & 0 & 0 \cr
                                           0 & \xi_2 & x_1 \cr
                                           0 & \ov{x}_1 & \xi_3},
                                  \pmatrix{\eta_1 & 0 & 0 \cr
                                           0 & -\eta_2 & -y_1 \cr
                                           0 & -\ov{y}_1 & -\eta_3}, -\xi, \eta \Big),
\vspace{1mm}\\
   \mu(X, Y, \xi, \eta) = \Big(\pmatrix{\eta & 0 & 0 \cr
                                        0 & \eta_3 & -y_1 \cr
                                        0 & -\ov{y}_1 & \eta_2},
                               \pmatrix{\xi & 0 & 0 \cr
                                        0 & \xi_3 & -x_1 \cr
                                        0 & -\ov{x}_1 & \xi_2}, 
        \eta_1, \xi_1 \Big).
\end{array}$$

We define subgroups $(E_7)^{\kappa,\mu}$, $((E_7)^{\kappa,\mu})_{(0,E_1,0,1)}$ and $((E_7)^{\kappa,\mu})_{(0,E_1,0,1),(0,-E_1,0,1)}$ of $E_7$ by
$$
      (E_7)^{\kappa,\mu} = \{\alpha \in E_7 \, | \, \kappa\alpha = \alpha\kappa, \mu\alpha = \alpha\mu \}, $$
\vspace{-3.5mm}
$$
    ((E_7)^{\kappa,\mu})_{(0,E_1,0,1)} = \{\alpha \in (E_7)^{\kappa,\mu} \, | \, \alpha(0, E_1, 0, 1) = (0, E_1, 0, 1) \}, $$
\vspace{-3mm}
$$
    ((E_7)^{\kappa,\mu})_{(0,E_1,0,1),(0,-E_1,0,1)} = \Big\{\alpha \in (E_7)^{\kappa,\mu} \; \Big| \begin{array}{l} \alpha(0, E_1, 0, 1) = (0, E_1, 0, 1)\\
                                   \alpha(0, -E_1, 0, 1) = (0, -E_1, 0, 1)
                  \end{array} \Big\}, $$
and also define subgroups $((E_7)^{\kappa,\mu})_{(E_1,0,1,0)}$ and $((E_7)^{\kappa,\mu})_{(E_1,0,1,0),(E_1,0,-1,0)}$ of $E_7$ by
$$
\begin{array}{c}
    ((E_7)^{\kappa,\mu})_{(E_1,0,1,0)} = \{\alpha \in (E_7)^{\kappa,\mu} \, | \, \alpha(E_1, 0, 1, 0) = (E_1, 0, 1, 0) \}, 
\vspace{1.5mm}\\
    ((E_7)^{\kappa,\mu})_{(E_1,0,1,0),(E_1,0,-1,0)} = \Big\{\alpha \in (E_7)^{\kappa,\mu} \; \Big| \begin{array}{l} \alpha(E_1, 0, 1, 0) = (E_1, 0, 1, 0)\\
                                   \alpha(E_1, 0, -1, 0) = (E_1, 0, -1, 0)
                  \end{array} 
\Big\}.
\vspace{2mm}
\end{array}$$

{\bf Proposition 3.2.} (1) $((E_7)^{\kappa,\mu})_{(E_1,0,1,0)} = ((E_7)^{\kappa,\mu})_{(0,E_1,0,1)}$.
\vspace{1mm}

(2) $((E_7)^{\kappa,\mu})_{(E_1,0,1,0),(E_1,0,-1,0)} = ((E_7)^{\kappa,\mu})_{(0,E_1,0,1),(0,-E_1,0,1)}$.
\vspace{2mm}

{\bf Proof.} (1) For $\alpha \in ((E_7)^{\kappa,\mu})_{(E_1,0,1,0)}$, we have $\alpha(0, E_1, 0, 1) = \alpha\mu(E_1, 0, 1, 0) = \mu\alpha(E_1, 0, 1, 0) = \mu(E_1, 0, 1, 0) = (0, E_1, 0, 1)$. Hence, $\alpha \in ((E_7)^{\kappa,\mu})_{(0,E_1,0,1)}$.
The converse is also proved.
\vspace{1mm}

(2) It is proved in a way similar to (1).
\vspace{3mm}

{\bf Proposition 3.3.} \qquad $((E_7)^{\kappa,\mu})_{(0,E_1,0,1),(0,-E_1,0,1)} \cong Spin(10)$.
\vspace{2mm}

{\bf Proof.} If $\alpha \in E_7$ satisfies $\alpha(0, E_1, 0, 1) = (0, E_1, 0, 1)$ and $\alpha(0, -E_1, 0, 1) = (0, -E_1, 0, 1)$, then we have $\alpha(0, 0, 0, 1) = (0, 0, 0, 1)$ and $\alpha(0, E_1, 0, 0) = (0, E_1, 0, 0)$. From the first condition, we see that $\alpha \in E_6$. Moreover, from the second condition, we have $\alpha \in (E_6)_{E_1} = Spin(10)$. The proof of the converse is trivial because $\kappa, \mu$ are defined by using $E_1$.   
\vspace{3mm}

{\bf Proposition 3.4.} \qquad \qquad $((E_7)^{\kappa,\mu})_{(0,E_1,0,1)} \cong Spin(11)$. 
\vspace{2mm}

{\bf Proof.} We define an 11 dimensional $\R$-vector space $V^{11}$ by
\begin{eqnarray*}
   V^{11} \!\!\! &=& \!\!\! \{P \in \gP^C \, | \, \kappa P = P, \mu\tau\lambda P = P, P \times (0, E_1, 0, 1) = 0 \}
\vspace{1mm}\\
        \!\!\! &=& \!\!\! \Big\{\Big(\pmatrix{0 & 0 & 0 \cr
                                              0 & \xi & x \cr
                                              0 & \ov{x} & -\tau\xi},
                                     \pmatrix{\eta & 0 & 0 \cr
                                              0 & 0 & 0 \cr
                                              0 & 0 & 0}, 0, \tau\eta \Big) 
   \Big| \, x \in \gC, \xi \in C, \eta \in i\R \Big\}
\end{eqnarray*}
with the norm 
$$
     (P, P)_\mu = \dfrac{1}{2}(\mu P, \lambda P) = (\tau\eta)\eta + \ov{x}x + (\tau\xi)\xi. $$
Let $SO(11) = SO(V^{11})$. Then, we have $((E_7)^{\kappa,\mu})_{(0,E_1,0,1)}/\Z_2 \cong SO(11), \Z_2 = \{1,\sigma \}$. Therefore, $((E_7)^{\kappa,\mu})_{(0,E_1,0,1)}$ is isomorphic to $Spin(11)$ as a double covering  group of $SO(11)$. (In details, see [6],[8].)
\vspace{2mm}

Now, we shall consider the following group
$$
\begin{array}{l}
    ((Spin(11))^{\sigma'})_{(0,F_1(y),0,0)}
\vspace{1mm}\\
\quad
   = \Big\{\alpha \in (Spin(11))^{\sigma'} \, \Big| \begin{array}{l} \alpha(0, F_1(y), 0, 0) \\
= (0, F_1(y), 0, 0) 
\end{array} \mbox{for all}\; y \in \gC
 \Big\}.
\end{array}$$
\vspace{-3mm}

{\bf Lemma 3.5.} {\it The Lie algebra $((\spin(11))^{\sigma'})_{(0,F_1(y),0,0)}$ of the group \\
$((Spin(11))^{\sigma'})_{(0,F_1(y),0,0)}$ is given by 
$$
\begin{array}{l}
   ((\spin(11))^{\sigma'})_{(0,F_1(y),0,0)}
\vspace{1mm}\\
\quad 
   = \Big\{{\mit\Phi}\Big(i\pmatrix{0 & 0 & 0 \cr
                                    0 & \epsilon & 0 \cr
                                    0 & 0 & -\epsilon}^{\!\!\sim},  
                           \pmatrix{0 & 0 & 0 \cr
                                    0 & \rho & 0 \cr
                                    0 & 0 & \tau\rho},
                     - \tau\pmatrix{0 & 0 & 0 \cr
                                    0 & \rho & 0 \cr
                                    0 & 0 & \tau\rho}, 0 \Big) \, \Big| \, 
   \epsilon \in \R, \rho \in C \Big\}.
\end{array}$$
In particular, we have}
$$
   \dim(((\spin(11))^{\sigma'})_{(0,F_1(y),0,0)}) 
= 3. $$
\vspace{-5mm}

{\bf Lemma 3.6.} {\it For $a \in \R$, the maps $\alpha_k(a) : \gP^C \to \gP^C, k = 1, 2, 3$ defined by
$$
   \alpha_k(a)\pmatrix{X \vspace{1mm}\cr
                      Y \vspace{1mm}\cr
                      \xi \vspace{1mm}\cr
                      \eta}
   = \pmatrix{(1 + (\cos a - 1)p_k)X - 2(\sin a)E_k \times Y + \eta(\sin a)E_k\vspace{1mm}\cr
   2(\sin a)E_k \times X + (1 + (\cos a - 1)p_k)Y - \xi(\sin a)E_k 
\vspace{1mm}\cr
   ((\sin a)E_k, Y) + (\cos a)\xi
\vspace{1mm}\cr
   (- (\sin a)E_k, X) + (\cos a)\eta} $$
belong to the group $E_7$, where $p_k : \gJ^C \to \gJ^C$ is defined by
$$
   p_k(X) = (X, E_k)E_k + 4E_k \times (E_k \times X), \quad X \in \gJ^C.$$
$\alpha_1(a), \alpha_2(b), \alpha_3(c)$ }($a, b, c \in \R$) {\it commute with each other.}
\vspace{2mm}

{\bf Proof.} For ${\mit\Phi}_k(a) = {\mit\Phi}(0, aE_k, -aE_k, 0) \in \gge_7$, we have $\alpha_k(a) = \exp{\mit\Phi}_k(a) \in E_7$. Since $[{\mit\Phi}_k(a), {\mit\Phi}_l(b)] = 0, k \not= l, \alpha_k(a)$ and $\alpha_l(b)$ are commutative.
\vspace{3mm}

{\bf Lemma 3.7.} \qquad \quad $((Spin(11))^{\sigma'})_{(0,F_1(y),0,0)}/Spin(2) \simeq S^2$.
\vspace{1mm}

\noindent {\it In particular, $((Spin(11))^{\sigma'})_{(0,F_1(y),0,0)}$ is connected.}
\vspace{2mm}

{\bf Proof.} We define a 3 dimensional $\R$-vector space $W^3$ by
\begin{eqnarray*}
   W^3 \!\!\! &=&\!\!\! \{P \in \gP^C \, |\, \kappa P = -P, \mu\tau\lambda P = -P, \sigma'P = P, P \times (E_1, 0, 1, 0) = 0 \}
\vspace{1mm}\\
    \!\!\! &=&\!\!\! \Big\{P = \Big(\pmatrix{i\xi & 0 & 0 \cr
                                             0 & 0 & 0 \cr
                                             0 & 0 & 0}, 
                                    \pmatrix{0 & 0 & 0 \cr
                                             0 & \eta & 0 \cr
                                             0 & 0 & -\tau\eta}, -i\xi, 0 \Big)
  \Big| \, \xi \in \R, \eta \in C \Big\}
\end{eqnarray*}
with the norm
$$
   (P, P)_\mu = - \dfrac{1}{2}(\mu P, \lambda P) = \xi^2 + (\tau\eta)\eta. $$
Then, $S^2 = \{P \in W^3 \, | \, (P, P)_\mu = 1 \}$ is a 2 dimensional sphere. The group \\
$((Spin(11))^{\sigma'})_{(0,F_1(y),0,0)}$ acts on $S^2$.  We shall show that this action is transitive. To show this, it is sufficient to show that any element $P \in S^2$ can be transformed to $(-iE_1, 0, i, 0) \in S^2$ under the action of $((Spin(11))^{\sigma'})_{(0,F_1(y),0,0)}$ . Now, for a given 
$$
   P = \Big(\pmatrix{i\xi & 0 & 0 \cr
                     0 & 0 & 0 \cr
                     0 & 0 & 0},
            \pmatrix{0 & 0 & 0 \cr
                     0 & \eta & 0 \cr
                     0 & 0 & -\tau\eta}, -i\xi, 0 \Big) \in S^2, $$
choose $a \in \R, 0 \le a < \pi/2$ such that $\tan2a = - \dfrac{2i\xi}{\tau\eta - \eta}$ (if $\tau\eta - \eta = 0$, then let $a 
\vspace{0.5mm}
= \pi/4)$. Operate $\alpha_{23}(a):= \alpha_2(a)\alpha_3(a) = \exp({\mit\Phi}(0, a(E_2 + E_3), -a(E_2 + E_3), 0))$ $\in ((Spin(11))^{\sigma'})_{(0,F_1(y),0,0)}$(Lemmas 3.5, 3.6) on $P$. Then, we have the $\xi$-term of $\alpha_{23}(a)P$ is $-((\cos2a)(i\xi) + 1/2 (\sin2a)(\tau\eta - \eta)) = 0$. Hence,
$$
   \alpha_{23}(a)P = \Big(0, \pmatrix{0 & 0 & 0 \cr
                                            0 & \zeta & 0 \cr
                                            0 & 0 & -\tau\zeta}, 0, 0 \Big) = P_1 , \quad  
      \zeta \in C,\,\, (\tau\zeta)\zeta = 1 . $$
From $(\tau\zeta)\zeta = 1, \zeta \in C$, we can put $\zeta = e^{i\theta}, 0 \le \theta < 2\pi$. Let $\nu = e^{-i\theta/2}$, and operate $\phi_1(\nu) \in ((Spin(10))^{\sigma'})_{F_1(x)}$ (Lemma 2.2) ($\subset ((Spin(11)^{\sigma'})_{(0,F_1(x),0,0)}$) on $P_1$. Then,
$$
    \phi_1(\nu)P_1 = (0, E_2 - E_3, 0, 0) = P_2. $$
Moreover, operate $\phi_1(e^{i\pi/4})$ on $P_2$,
$$  
   \phi_1(e^{i\pi/4})P_2 = (0, i(E_2 + E_3), 0, 0) = P_3. $$
Operate again $\alpha_{23}(\pi/4)$ on $P_3$. Then, we have
$$
    \alpha_{23}(\pi/4)P_3 = (-iE_1, 0, i, 0). $$
This shows the transitivity. The isotropy subgroup of $((Spin(11))^{\sigma'})_{(0,F_1(y),0,0)}$ at $(-iE_1, 0, i, 0)$ is $((Spin(10))^{\sigma'})_{F_1(y)}$ (Propositions 3.2(2), 3.3, 3.4) $ = Spin(2)$. Thus, we have the homeomorphism $((Spin(11))^{\sigma'})_{(0,F_1(y),0,0)}/Spin(2) \simeq S^2$. 
\vspace{3mm}

{\bf Proposition 3.8.} \qquad \quad $((Spin(11))^{\sigma'})_{(0,F_1(y),0,0)} \cong Spin(3)$.
\vspace{2mm}

{\bf Proof.}  Since $((Spin(11))^{\sigma'})_{(0,F_1(y),0,0)}$ is connected (Lemma 3.7), we can define a homomorphism $\pi : ((Spin(11))^{\sigma'})_{(0,F_1(y),0,0)} \to SO(3) = SO(W^3)$ by
$$
                \pi(\alpha) = \alpha|W^3. $$
$\Ker\,\pi = \{1, \sigma \} = \Z_2$. \, Since $\dim(((\spin(11))^{\sigma'})_{(0,F_1(y),0,0)}) = \,3 \,\mbox{(Lemma 3.5)}\, =$ \\  
$\dim(\so(3)),\, \pi$ is onto. \, Hence, $((Spin(11))^{\sigma'})_{(0,F_1(y),0,0)}/\Z_2 \cong SO(3)$. \,Therefore, \\
$((Spin(11))^{\sigma'})_{(0,F_1(y),0,0)}$ is isomorphic to $Spin(3)$ as a double covering group of $SO(3)$. 
\vspace{3mm}

{\bf Lemma 3.9.} {\it The Lie algebra $(\spin(11))^{\sigma'}$ of the group $(Spin(11))^{\sigma'}$ is given by} $$
\begin{array}{l}
   (\spin(11))^{\sigma'}
\vspace{1mm}\\
\qquad
   = \Big\{{\mit\Phi}\Big(D + i\pmatrix{0 & 0 & 0 \cr
                                      0 & \epsilon & 0 \cr
                                      0 & 0 & -\epsilon}^{\!\!\sim},
                             \pmatrix{0 & 0 & 0 \cr
                                      0 & \rho & 0 \cr
                                      0 & 0 & \tau\rho},
                        -\tau\pmatrix{0 & 0 & 0 \cr
                                      0 & \rho & 0 \cr
                                      0 & 0 & \tau\rho}, 0 \Big)
\vspace{1mm}\\
\qquad \quad
   \Big| \, D \in \so(8), \epsilon \in \R, \rho \in C \Big\}.
\end{array}$$
{\it In particular, we have}
$$
           \dim((\spin(11))^{\sigma'}) = 
28 + 3 = 31. $$
\vspace{-5mm} 

Now, we shall determine the group structure of $(Spin(11))^{\sigma'}$.
\vspace{3mm}

{\bf Theorem 3.10.} \enskip $(Spin(11))^{\sigma'} \cong (Spin(3) \times Spin(8))/\Z_2, \Z_2 = \{(1, 1), (-1, \sigma)\}.$
\vspace{-3mm}

{\bf Proof.} Let $Spin(11) = ((E_7)^{\kappa,\mu})_{(0,E_1,0,1)}$, $Spin(3) = ((Spin(11))^{\sigma'})_{(0,F_1(y),0,0)}$ and $Spin(8) = ((F_4)_{E_1})^{\sigma'} \subset ((E_6)_{E_1})^{\sigma'} = (((E_7)^{\kappa,\mu})_{(E_1,0,1,0),(E_1,0,-1,0)})^{\sigma'} \subset \\
(((E_7)^{\kappa,\mu})_{(E_1,0, 1, 0)})^{\sigma'}$(Theorem 1.2, Propositions 3.2, 3.3, 3.4). Now, we define a map $\varphi : Spin(3) \times Spin(8) \to (Spin(11))^{\sigma'}$ by
$$
   \varphi(\alpha, \beta) = \alpha\beta. $$
Then, $\varphi$ is well-defined : $\varphi(\alpha, \beta) \in (Spin(11))^{\sigma'}$. Since
\vspace{1mm}
 $[{\mit\Phi}_D, {\mit\Phi}_3] = 0$ for ${\mit\Phi}_D = {\mit\Phi}(D, 0, 0, 0) \in \spin(8)$, ${\mit\Phi}_3 \in \spin(3) = ((\spin(11))^{\sigma'})_{(0,F_1(y),0,0)}$ (Proposition  3.8), we have $\alpha\beta = \beta\alpha$. Hence, $\varphi$ is a homomorphism. $\Ker\,\varphi = \{(1, 1), (-1, \sigma) \} = \Z_2$. Since $(Spin(11))^{\sigma'}$ is connected and $\dim(\spin(3) \oplus \spin(8)) = 3 \mbox{(Lemma 3.5)} + 28 = 31 = \dim((\spin(11))^{\sigma'})$ (Lemma 3.9), $\varphi$ is onto. Thus, we have the isomorphism $(Spin(3) \times Spin(8))/\Z_2 \cong (Spin(11))^{\sigma'} $.
\vspace{3mm}

{\bf Proposition 3.11.} \qquad \qquad $(E_7)^{\kappa,\mu} \cong Spin(12)$.
\vspace{2mm}

{\bf Proof.} We define a 12 dimensional $\R$-vector space $V^{12}$ by
\begin{eqnarray*}
   V^{12} \!\!\! &=& \!\!\! \{P \in \gP^C \, | \, \kappa P = P, \mu\tau\lambda P = P \}
\vspace{1mm}\\
        \!\!\! &=& \!\!\! \Big\{\Big(\pmatrix{0 & 0 & 0 \cr
                                              0 & \xi & x \cr
                                              0 & \ov{x} & -\tau\xi},
                      \pmatrix{\eta & 0 & 0 \cr
                               0 & 0 & 0 \cr
                               0 & 0 & 0}, 0, \tau\eta \Big) \Big| \, x \in \gC, \xi, \eta \in C \Big\}
\end{eqnarray*}
with the norm 
$$
    (P, P)_\mu = \dfrac{1}{2}(\mu P, \lambda P) = (\tau\eta)\eta + \ov{x}x + (\tau\xi)\xi. $$
Let $SO(12) = SO(V^{12})$. Then, we have $(E_7)^{\kappa,\mu}/ \Z_2 \cong SO(12)$, $\Z_2 = \{1, \sigma \}$. Therefore, $(E_7)^{\kappa,\mu}$ is isomorphic to $Spin(12)$ as a double covering group of $SO(12)$. (In details, see [6],[8].)
\vspace{2mm}

Now, we shall consider the following group
$$
\begin{array}{l}
    ((Spin(12))^{\sigma'})_{(0,F_1(y),0,0)}
\vspace{1mm}\\
\quad
   = \Big\{\alpha \in (Spin(12))^{\sigma'} \, \Big| \begin{array}{l} \alpha(0, F_1(y), 0, 0) \\
= (0, F_1(y), 0, 0) 
\end{array} \mbox{for all}\; y \in \gC
 \Big\}.
\end{array}$$
\vspace{-3mm}

{\bf Lemma 3.12.} {\it The Lie algebra $((\spin(12))^{\sigma'})_{(0,F_1(y),0,0)}$ of the group}\\  
$((Spin(12))^{\sigma'})_{(0,F_1(y),0,0)}$ {\it is given by }
$$
\begin{array}{l}
   ((\spin(12))^{\sigma'})_{(0,F_1(y),0,0)}
\vspace{1mm}\\
\qquad
   = \Big\{{\mit\Phi}\Big(i\pmatrix{\epsilon_1 & 0 & 0 \cr
                                    0 & \epsilon_2 & 0 \cr
                                    0 & 0 & \epsilon_3}^{\!\!\sim},  
                           \pmatrix{0 & 0 & 0 \cr
                                    0 & \rho_2 & 0 \cr
                                    0 & 0 & \rho_3},
                     - \tau\pmatrix{0 & 0 & 0 \cr
                                    0 & \rho_2 & 0 \cr
                                    0 & 0 & \rho_3}, - \dfrac{3}{2}i\epsilon_1 \Big) 
\vspace{1mm}\\
\qquad \quad
 \Big| \, \epsilon_i \in \R, \epsilon_1 + \epsilon_2 + \epsilon_3 = 0, \rho_i \in C \Big\}.
\end{array}$$
{\it In particular, we have}
$$
        \dim(((\spin(12))^{\sigma'})_{(0,F_1(y),0,0)}) = 6. $$

{\bf Lemma 3.13.} {\it For $t \in \R$, the map $\alpha(t) : \gP^C \to \gP^C$ defined by
$$
\begin{array}{l}
   \alpha(t)(X, Y, \xi, \eta)
\vspace{1mm}\\
\quad
   = \Big(\pmatrix{e^{2it}\xi_1 & e^{it}x_3 & e^{it}\ov{x}_2 \cr
                    e^{it}\ov{x}_3 & \xi_2 & x_1 \cr
                    e^{it}x_2 & \ov{x}_1 & \xi_3}, 
           \pmatrix{e^{-2it}\eta_1 & e^{-it}y_3 & e^{-it}\ov{y}_2 \cr
                    e^{-it}\ov{y}_3 & \eta_2 & y_1 \cr
                    e^{-it}y_2 & \ov{y}_1 & \eta_3}, e^{-2it}\xi, e^{2it}\eta \Big)
\end{array}$$
belongs to the group $((Spin(12))^{\sigma'})_{(0,F_1(y),0,0)}$.}
\vspace{2mm}

{\bf Proof.} For ${\mit\Phi} = {\mit\Phi}(2itE_1 \vee E_1, 0, 0, -2it) \in ((\spin(12))^{\sigma'})_{(0,F_1(y),0,0)}$ (Lemma 3.12), we have $\alpha(t) = \exp{\mit\Phi} \in ((Spin(12)^{\sigma'})_{(0,F_1(y),0,0)}$. 
\vspace{3mm}

{\bf Lemma 3.14.} \qquad \quad $((Spin(12))^{\sigma'})_{(0,F_1(y),0,0)}/Spin(3) \simeq S^3$.
\vspace{1mm}

\noindent {\it In particular, $((Spin(12))^{\sigma'})_{(0,F_1(y),0,0)}$ is connected.}
\vspace{2mm}

{\bf Proof.} We define a 4 dimensional $\R$-vector space $W^4$ by
\begin{eqnarray*}
   W^4 \!\!\! &=&\!\!\! \{P \in \gP^C \, |\, \kappa P = -P, \mu\tau\lambda P = -P, \sigma'P = P \}
\vspace{1mm}\\
    \!\!\! &=& \!\!\! \Big\{P = \Big(\pmatrix{\xi & 0 & 0 \cr
                                           0 & 0 & 0 \cr
                                           0 & 0 & 0}, 
                                  \pmatrix{0 & 0 & 0 \cr
                                           0 & \eta & 0 \cr
                                           0 & 0 & -\tau\eta}, \tau\xi, 0 \Big)
  \Big| \, \xi, \eta \in C \Big\}
\end{eqnarray*}
with the norm
$$
   (P, P)_\mu = - \dfrac{1}{2}(\mu P, \lambda P) = (\tau\xi)\xi + (\tau\eta)\eta. $$
Then, $S^3 = \{P \in W^4 \, | \, (P, P)_\mu = 1 \}$ is a 3 dimensional sphere. The group \\
$((Spin(12))^{\sigma'})_{(0,F_1(y),0,0)}$ acts on $S^3$.  We shall show that this action is transitive. To show this, it is sufficient to show that any element $P \in S^3$ can be transformed to $(E_1, 0, 1, 0) \in S^3$ under the action of $((Spin(12))^{\sigma'})_{(0,F_1(y),0,0)}$. Now, for a given 
$$
   P = \Big(\pmatrix{\xi & 0 & 0 \cr
                     0 & 0 & 0 \cr
                     0 & 0 & 0},
            \pmatrix{0 & 0 & 0 \cr
                     0 & \eta & 0 \cr
                     0 & 0 & -\tau\eta}, \tau\xi, 0 \Big) \in S^3, $$
choose $t \in \R$ such that $e^{2it}\xi \in i\R$. Operate $\alpha(t)$ (Lemma 3.13) on $P$. Then, we have
$$
            \alpha(t)P = P_1 \in S^2 \subset S^3 . $$
Now, since $((Spin(11))^{\sigma'})_{(0,F_1(y),0,0)}$ ($\subset ((Spin(12))^{\sigma'})_{(0,F_1(y),0,0)}$) acts transitively on $S^2$ (Lemma 3.7), there exists $\beta \in ((Spin(11))^{\sigma'})_{(0,F_1(y),0,0)}$ such that 
$$
         \beta P_1 = (-iE_1, 0, i, 0) = P_2. $$
Operate again $\alpha(\pi/4)$ on $P_2$. Then, we have
$$
       \alpha(\pi/4)P_2 = (E_1, 0, 1, 0). $$
This shows the transitivity. The isotropy subgroup of $((Spin(12))^{\sigma'})_{(0,F_1(y),0,0)}$ at $(E_1, 0, 1, 0)$ is $((Spin(11))^{\sigma'})_{(0,F_1(y),0,0)}$ (Propositions 3.2(1), 3.4, 3.11) $ = Spin(3)$. Thus, we have the homeomorphism $((Spin(12))^{\sigma'})_{(0,F_1(y),0,0)}/Spin(3) \simeq S^3$. 
\vspace{3mm}

{\bf Proposition 3.15.} \qquad \quad $((Spin(12))^{\sigma'})_{(0,F_1(y),0,0)} \cong Spin(4)$.
\vspace{2mm}

{\bf Proof.} Since $((Spin(12))^{\sigma'})_{(0,F_1(y),0,0)}$ is connected (Lemma 3.14), we can define a homomorphism $\pi : ((Spin(12))^{\sigma'})_{(0,F_1(y),0,0)} \to SO(4) = SO(W^4)$ by
$$
                \pi(\alpha) = \alpha|W^4. $$
$\Ker\,\pi\!=\! \{1, \sigma \} \!=\! \Z_2$. Since $\dim((\spin(12))^{\sigma'})_{(0,F_1(y),0,0)}) = 6$ (Lemma 3.12) $= \dim(\so(4))$, $\pi$ is onto. Hence, $((Spin(12))^{\sigma'})_{(0,F_1(y),0,0)}/\Z_2 \cong SO(4)$. Therefore, $((Spin(12))^{\sigma'})_{(0,F_1(y),0,0)}$ is isomorphic to $Spin(4)$ as a double covering group of 
\vspace{3mm}
$SO(4)$. 

{\bf Lemma 3.16.} {\it The Lie algebra $(\spin(12))^{\sigma'}$ of the group $(Spin(12))^{\sigma'}$ is given by} 
$$
\begin{array}{l}
   (\spin(12))^{\sigma'}
\vspace{1mm}\\
\qquad
   = \Big\{{\mit\Phi}\Big(D + i\pmatrix{\epsilon_1 & 0 & 0 \cr
                                        0 & \epsilon_2 & 0 \cr
                                        0 & 0 & \epsilon_3}^{\!\!\sim},
                               \pmatrix{0 & 0 & 0 \cr
                                        0 & \rho_2 & 0 \cr
                                        0 & 0 & \rho_3},
                          -\tau\pmatrix{0 & 0 & 0 \cr
                                        0 & \rho_2 & 0 \cr
                                        0 & 0 & \rho_3}, -i\dfrac{3}{2}\epsilon_1 \Big)
\vspace{1mm}\\
\qquad \quad
   \Big| \, D \in \so(8), \epsilon_i \in \R, \epsilon_1 + \epsilon_2 + \epsilon_3 = 0, \rho_i \in C \Big\}.
\end{array}$$
{\it In particular, we have}
$$
           \dim((\spin(12))^{\sigma'}) = 
\vspace{2mm}
28 + 6 = 34. $$ 

Now, we shall determine the group structure of $(Spin(12))^{\sigma'}$.
\vspace{3mm}

{\bf Theorem 3.17.} \enskip $(Spin(12))^{\sigma'} \cong (Spin(4) \times Spin(8))/\Z_2, \Z_2 = \{(1, 1), (-1, \sigma) \}.$
\vspace{-3mm}

{\bf Proof.} Let $Spin(12) = (E_7)^{\kappa,\mu}$, $Spin(4) = ((Spin(12))^{\sigma'})_{(0,F_1(y),0,0)}$ and $Spin(8)$  $ = ((F_4)_{E_1})^{\sigma'} \subset ((E_6)_{E_1})^{\sigma'} = (((E_7)^{\kappa,\mu})_{(E_1,0,1,0),(E_1,0,-1,0)})^{\sigma'} \subset ((E_7)^{\kappa,\mu})^{\sigma'}$(Theorem 1.2, Propositions 3.2, 3.3, 3.11, 3.15). Now, we define a map $\varphi : Spin(4) \times Spin(8) \to (Spin(12))^{\sigma'}$ by
$$
         \varphi(\alpha, \beta) = \alpha\beta. $$
Then, $\varphi$ is well-defined : $\varphi(\alpha, \beta) \in (Spin(12))^{\sigma'}$. Since
 $[{\mit\Phi}_D, {\mit\Phi}_4] = 0$ for ${\mit\Phi}_D = {\mit\Phi}(D, 0, 0, 0) \in \spin(8)$, ${\mit\Phi}_4  \in
\vspace{1mm}
 \spin(4) = ((\spin(12))^{\sigma'})_{(0,F_1(y),0,0)}$ (Proposition  3.15), we have $\alpha\beta = \beta\alpha$. Hence, $\varphi$ is a homomorphism. $\Ker\,\varphi = \{(1, 1), (-1, \sigma) \} = \Z_2$. Since $(Spin(12))^{\sigma'}$ is connected and $\dim(\spin(4) \oplus \spin(8)) = 6 \mbox{(Lemma 3.12)} + 28 = 34 = \dim((\spin(12))^{\sigma'})$ (Lemma 3.16), $\varphi$ is onto. Thus, we have the isomorphism $(Spin(4) \times Spin(8))/\Z_2 \cong (Spin(12))^{\sigma'}$.
\vspace{4mm}

\begin{center}
                  {\large{\bf 4. Group $E_8$}}
\end{center}
\vspace{3mm}

We use the same notation as in [2],[4](however, some will rewritten). For example,
\vspace{-4mm}

$C$-Lie algebra ${\gge_8}^C = {\gge_7}^C \oplus \gP^C \oplus \gP^C \oplus C \oplus C \oplus C$ and $C$-linear transformations $\lambda, \wti{\lambda}$ of ${\gge_8}^C$,
\vspace{1mm}

the groups ${E_8}^C = \{\alpha \in \Iso_C({\gge_8}^C) \, | \, \alpha[R_1, R_2] = [\alpha R_1, \alpha R_2] \}$ and $E_8 = ({E_8}^C)^{\tau\wti{\lambda}} = \{\alpha \in {E_8}^C \, | \, \tau\wti{\lambda}\alpha = \alpha\tau\wti{\lambda} \}$.
\vspace{3mm}

For $\alpha \in E_7$, the map $\wti{\alpha} : {\gge_8}^C \to {\gge_8}^C$ is defined by
$$ 
     \wti{\alpha}({\mit\Phi}, P, Q, r, u, v) = (\alpha{\mit\Phi}{\alpha}^{-1}, \alpha P, \alpha Q, r, u, v).$$
Then, $\wti{\alpha} \in E_8$ and we identify $\alpha$ with $\wti{\alpha}$. The group $E_8$ contains $E_7$ as a subgroup by 
$$    
   E_7 = \{\wti{\alpha} \in E_8 \,| \, \alpha \in E_7 \} = (E_8)_{(0, 0, 0, 0, 1, 0)}.$$

We define a $C$-linear map $\wti{\kappa} : {\gge_8}^C \to {\gge_8}^C$ by
$$
   \wti{\kappa} = \mbox{ad}(\kappa, 0, 0, -1, 0, 0 ) = \mbox{ad}({\mit\Phi}(-2 E_1 \vee E_1, 0, 0, -1), 0 ,0, -1, 0, 0 ), $$ 
and 14 dimensional $C$-vector spaces $\gg_{-2}$ and $\gg_2$ by
\begin{eqnarray*}
   \gg_{-2} \!\!\! &=& \!\!\! \{R \in {\gge_8}^C\,|\, \wti{\kappa}R = -2R \},
\vspace{1mm}\\
     \!\!\! &=& \!\!\! \{({\mit\Phi}(0, \zeta E_1, 0, 0), (\xi_1E_1, \eta_2E_2 + \eta_3E_3 + F_1(y), \xi, 0), 0, 0, u, 0)
\vspace{1mm}\\
   & & | \, \zeta, \xi_1, \eta_i, \xi, u \in C, y \in \gC^C \},
\vspace{1mm}\\
   \gg_2 \!\!\! &=& \!\!\! \{R \in {\gge_8}^C\,|\, \wti{\kappa}R = 2R \},
\vspace{1mm}\\
     \!\!\! &=& \!\!\! \{({\mit\Phi}(0, 0, \zeta E_1, 0), 0, (\xi_2E_1 + \xi_3E_3 + F_1(x), \eta_1E_1, 0, \eta), 0, 0, v)
\vspace{1mm}\\
   & & | \, \zeta, \xi_i, \eta_1, \eta, v \in C, x \in \gC^C \}.
\end{eqnarray*}
Further, we define two $C$-linear maps $\wti{\mu}_1 : {\gge_8}^C \to {\gge_8}^C$ and $\delta : \gg_2 \to \gg_2$ by 
$$ 
   \wti{\mu}_1({\mit\Phi}, P, Q, r, u, v) = (\mu_1{\mit\Phi}{\mu_1}^{-1}, i\mu_1 Q, i\mu_1 P, -r, v, u),$$
\noindent where 
$$     
    \mu_1(X, Y, \xi, \eta) =
       \Bigl(\pmatrix{i\eta & x_3 & \ov{x}_2 \vspace{0.5mm}\cr
                      \ov{x}_3 & i\eta_3 & - iy_1 \vspace{0.5mm}\cr
                      x_2 & - i\ov{y}_1 & i\eta_2},
             \pmatrix{i\xi & y_3 & \ov{y}_2 \vspace{0.5mm}\cr
                      \ov{y}_3 & i\xi_3 & - ix_1 \vspace{0.5mm}\cr
                      y_2 & - i\ov{x}_1 & i\xi_2},
                 i\eta_1, i\xi_1 \Bigr), $$
\noindent and 
$$
\begin{array}{l} 
    \delta({\mit\Phi}(0, 0, \zeta E_1, 0), 0, (\xi_2E_2 + \xi_3E_3 + F_1(x), \eta_1E_1, 0, \eta), 0, 0, v) 
\vspace{1mm}\\
\quad
   = ({\mit\Phi}(0, 0, -vE_1, 0), 0, (\xi_2E_2 + \xi_3E_3 + F_1(x), \eta_1E_1, 0, \eta), 0, 0, -\zeta). 
\end{array}$$
In particular, the explicit form of the map $\wti{\mu}_1 : \gg_{-2} \to \gg_2$ is given by
$$
\begin{array}{l}
    \wti{\mu}_1({\mit\Phi}(0, \zeta E_1, 0, 0), (\xi_1E_1, \eta_2E_2 + \eta_3E_3 + F_1(y), \xi, 0), 0, 0, u, 0) 
\vspace{1mm}\\
\quad
   = ({\mit\Phi}(0, 0, \zeta E_1, 0), 0, (-\eta_3E_2 - \eta_2E_3 + F_1(y), -\xi E_1, 0, -\xi_1), 0, 0, u). 
\end{array}$$
The composition map $\delta\wti{\mu}_1 : \gg_{-2} \to \gg_2$ of $\wti{\mu}_1$ and $\delta\wti{\mu}_1$ is denoted by $\wti{\mu}_{\delta}$:
$$
\begin{array}{l} 
  \wti{\mu}_{\delta}({\mit\Phi}(0, \zeta E_1, 0, 0), (\xi_1E_1, \eta_2E_2 + \eta_3E_3 + F_1(y), \xi, 0), 0, 0, u, 0)
\vspace{1mm}\\
\quad
   = ({\mit\Phi}(0, 0, -uE_1, 0), 0, (-\eta_3E_2 - \eta_2E_3 + F_1(y), -\xi E_1, 0, -\xi_1), 0, 0, -\zeta). 
\end{array}$$
Now, we define the inner product $(R_1, R_2)_{\mu}$ in $\gg_{-2}$ by
$$ 
        (R_1, R_2)_{\mu} = \frac{1}{30}B_8(\wti{\mu}_{\delta}R_1, R_2), $$
where $B_8$ is the Killing form of ${\gge_8}^C$. The explicit form of $(R, R)_{\mu}$ is given by
$$
   (R, R)_{\mu} = -4\zeta u - \eta_2\eta_3 + \ov{y}y + \xi_1\xi, $$
for $R = ({\mit\Phi}(0, \zeta E_1, 0, 0), (\xi_1 E_1, \eta_2E_2 + \eta_3E_3 + F_1(y), \xi, 0), 0, 0, u, 0) \in \gg_{-2}$.
\vspace{1mm}
Hereafter, we use the notation $(V^C)^{14}$ instead of $\gg_{-2}$.
\vspace{2mm}
              
We define $\R$-vector spaces $V^{14}, V^{13}$ and $(V')^{12}$ respectively by
\begin{eqnarray*}
     V^{14} \!\!\! &=& \!\!\! \{ R \in (V^C)^{14} \,|\, \wti{\mu}_{\delta}\tau\wti{\lambda}R = -R \}
\vspace{1mm}\\
     \!\!\! &=& \!\!\! \{R = ({\mit\Phi}(0, \zeta E_1, 0, 0), (\xi E_1, \eta E_2 - \tau\eta E_3 + F_1(y), \tau\xi, 0), 0, 0, -\tau\zeta, 0)
\vspace{1mm}\\
    & & \quad |\, \zeta, \xi, \eta \in C, y \in \gC \}
\end{eqnarray*}
with the norm
$$
   (R, R)_\mu = \dfrac{1}{30}B_8(\wti{\mu}_\delta R, R) = 4(\tau\zeta)\zeta + (\tau\eta)\eta + \ov{y}y + (\tau\xi)\xi, $$
\vspace{-4mm}
\begin{eqnarray*}
    V^{13} \!\!\! &=& \!\!\! \{ R \in V^{14} \,| \, (R, ({\mit\Phi}_1, 0, 0, 0, 1, 0))_{\mu} = 0 \}
\vspace{1mm}\\   
     \!\!\! &=& \!\!\! \{R = ({\mit\Phi}(0, \zeta E_1, 0, 0), (\xi E_1, \eta E_2 - \tau\eta E_3 + F_1(y), \tau\xi, 0), 0, 0, -\zeta, 0)
\vspace{1mm}\\
     & & \quad | \, \zeta \in \R, \xi, \eta \in C, y \in \gC \}
\end{eqnarray*}
with the norm
$$
   (R, R)_\mu = \dfrac{1}{30}B_8(\wti{\mu}_\delta R, R) = 4\zeta^2 + (\tau\eta)\eta + \ov{y}y + (\tau\xi)\xi, $$
\vspace{-4mm}
\begin{eqnarray*}
    (V')^{12} \!\!\! &=& \!\!\! \{ R \in V^{13} \,| \, (R, ({\mit\Phi}_1, 0, 0, 0, -1, 0))_{\mu} = 0 \}
\vspace{1mm}\\
   \!\!\! &=& \!\!\! \{R = (0, (\xi E_1, \eta E_2 - \tau\eta E_3 + F_1(y), \tau\xi, 0), 0, 0, 0, 0) 
\vspace{1mm}\\
     & & \quad | \, \xi, \eta \in C, y \in \gC \}
\end{eqnarray*}
with the norm
$$
   (R, R)_\mu = \dfrac{1}{30}B_8(\wti{\mu}_\delta R, R) = (\tau\eta)\eta + \ov{y}y + (\tau\xi)\xi, $$
where ${\mit\Phi}_1 = {\mit\Phi}(0 ,E_1 ,0, 0 )$. We use the notation $(V')^{12}$ to distinguish from the $\R$-vector space $V^{12}$ defined in Section 3. The space $(V')^{12}$ above can be identified with the $\R$-vector space
$$
\begin{array}{l}
   \quad  \{ P \in \gP^C \,| \, \kappa P = -P, \mu\tau\lambda P = -P \}
\vspace{1mm}\\

   = \{P = (\xi E_1, \eta E_2 - \tau\eta E_3 + F_1(y), \tau\xi, 0) \in \gP^C \,| \, \xi, \eta \in C, y \in \gC \}
\end{array}$$
with the norm
$
   (P, P)_\mu = - \dfrac{1}{2}(\mu P, \lambda P) = (\tau\eta)\eta + \ov{y}y + (\tau\xi)\xi. $
\vspace{2mm}

Now, we define a subgroup $G_{14}$ of ${E_8}^C$ by
$$ 
     G_{14} = \{ \alpha \in {E_8}^C \,|\, \wti{\kappa}\alpha = \alpha\wti{\kappa},\,\wti{\mu}_{\delta}
\alpha R = \alpha\wti{\mu}_{\delta}R,\, R \in (V^C)^{14} 
\vspace{1mm}
\}. $$

{\bf Lemma 4.1.} {\it The Lie algebra $\gg_{14}$ of the group $G_{14}$ is given by
$$
\begin{array}{l}
     \gg_{14} = \{R \in {\gge_8}^C \, | \, \wti{\kappa}(\ad R) = (\ad R) \wti{\kappa},\,(\wti{\mu}_{\delta}(\ad R))R'                                    =((\ad R)\wti{\mu}_{\delta})R',\,R' \in (V^C)^{14}\} 
\vspace{1mm}\\
\qquad
  = \Big\{\Bigl({\mit\Phi}\Bigl(D + \pmatrix{0 & 0 & 0 \cr
                                0 & 0 & d_1 \cr
                                0 & - \ov{d}_1 & 0}^{\!\!\sim} + 
                       \pmatrix{\tau_1 & 0 & 0 \cr
                                0 & \tau_2 & t_1 \cr
                                0 & \ov{t}_1 & \tau_3}^{\!\!\sim},
                       \pmatrix{0 & 0 & 0 \cr
                                0 & \alpha_2 & a_1 \cr
                                0 & \ov{a}_1 & \alpha_3},
\vspace{1mm}\\
\qquad
                       \pmatrix{0 & 0 & 0 \cr
                                0 & \beta_2 & b_1 \cr
                                0 & \ov{b}_1 & \beta_3}, \nu \Bigl), 

     \Bigl(\pmatrix{0 & 0 & 0 \cr
                   0 & \rho_2 & p_1 \cr
                   0 & \ov{p}_1 & \rho_3},
         \pmatrix{\rho_1 & 0 & 0 \cr
                   0 & 0 & 0 \cr
                   0 & 0 & 0}, 0, \rho \Bigl),
          \Bigl(\pmatrix{\zeta_1 & 0 & 0 \cr
                   0 & 0 & 0 \cr
                   0 & 0 & 0},
\vspace{1mm}\\
\qquad
          \pmatrix{0 & 0 & 0 \cr
                   0 & \zeta_2 & z_1 \cr
                   0 & \ov{z}_1 & \zeta_3},
                   \zeta, 0 \Bigl), r, 0, 0 \Bigl) 

 \Big| \, D \in \so(8)^C,\, \tau_i, \alpha_i, \beta_i, \nu, \rho_i, \rho, \zeta_i, \zeta, r \in C,
\vspace{1mm}\\
\qquad
    \tau_1 + \tau_2 + \tau_3 = 0, \, d_1, t_1, a_1, b_1, p_1, z_1 \in \gC^C, \,  \tau_1 + \dfrac{2}{3}\nu + 2r = 0 \Big\} . 
\end{array}$$
In particular, we have}
$$
    \dim_C(\gg_{14}) = 
\vspace{3mm}
28 + 63 = 91. $$

{\bf Proposition 4.2.}  \qquad \qquad $G_{14} \cong Spin(14, C)$.
\vspace{3mm}

{\bf Proof.} Let $SO(14, C) = SO((V^{14})^C)$. Then, we have $G_{14}/\Z_2 \cong SO(14, C), \Z_2 = \{1, \sigma \}$. Therefore, $G_{14}$ is isomorphic to $Spin(14, C)$ as a double covering group of $SO(14, C)$. (In details, see [2].) 
\vspace{2mm}

We define subgroups ${G_{14}}^{com}, {G_{13}}^{com}$ and ${G_{12}}^{com}$ of the group $E_8$ by 
\begin{eqnarray*}
 {G_{14}}^{com} \!\!\! &=& \!\!\! \{\alpha \in G_{14} \, | \, \tau\wti{\lambda}\alpha = \alpha\tau \wti{\lambda}\},
\vspace{1mm}\\
 {G_{13}}^{com} \!\!\! &=& \!\!\! \{\alpha \in {G_{14}}^{com}\,|\,\alpha({\mit\Phi}_1, 0, 0, 0, 1, 0) = ({\mit\Phi}_1, 0, 0, 0, 1, 0) \},
\vspace{1mm}\\
 {G_{12}}^{com} \!\!\! &=& \!\!\! \{\alpha \in {G_{13}}^{com}\,|\,\alpha({\mit\Phi}_1, 0, 0, 0, -1, 0) = ({\mit \Phi}_1, 0, 0, 0, -1, 0)\},
\end{eqnarray*}
\noindent respectively.
\vspace{3mm}

{\bf Lemma 4.3.} $\alpha \in (E_7)^{\kappa,\mu} = Spin(12)$ {\it satisfies} 
$$ 
       \alpha {\mit\Phi}(0, E_1, 0, 0 )\alpha^{-1} = {\mit\Phi}(0, E_1, 0, 0 ),
\quad and \quad
       \alpha {\mit\Phi}(0, 0, E_1, 0 )\alpha^{-1} = {\mit\Phi}(0, 0, E_1, 0 ).  $$

{\bf Proof}. We consider an 11 dimensional sphere $(S')^{11}$ by
$$
\begin{array}{l}
  (S')^{11} = \{ P' \in (V')^{12} \,| \, (P, P)_{\mu} = 1 \}
\vspace{1mm}\\
\qquad \quad \!
   = \{ P' = (\xi E_1, \eta E_2 - \tau\eta E_3 + F_1(y), \tau\xi, 0) 
\vspace{1mm}\\
\qquad \qquad \qquad \quad
     \,| \, \xi, \eta \in C,\, y \in \gC, (\tau\eta)\eta + \ov{y}y + (\tau\xi)\xi = 1 \}.
\end{array}$$
Since the group $Spin(12)$ acts on $(S')^{11}$, we can put 
$$ 
   \alpha(E_1, 0, 1, 0) = (\xi E_1, \eta E_2 - \tau\eta E_3 + F_1(y), \tau\xi, 0) \ \in (S')^{11}.$$
Now, since $1/2{\mit\Phi}(0, E_1, 0, 0) = (E_1, 0, 1, 0) \times (E_1, 0, 1, 0)$, we have 
$$
\begin{array}{l}
     1/2\alpha{\mit\Phi}(0, E_1, 0, 0)\alpha^{-1} = 
\alpha((E_1, 0, 1, 0) \times (E_1, 0, 1, 0))\alpha^{-1} 
\vspace{1.5mm}\\
\qquad
     = \alpha(E_1, 0, 1, 0) \times \alpha (E_1, 0, 1, 0)
\vspace{1.5mm}\\
\qquad
   = (\xi E_1, \eta E_2 - \tau\eta E_3 + F_1(y), \tau\xi, 0) \times (\xi E_1, \eta E_2 - \tau\eta E_3 + F_1(y), \tau\xi, 0) 
\vspace{1.5mm}\\
\qquad 
   = 1/2{\mit\Phi}(0, ((\tau\eta)\eta + \ov{y}y + (\tau\xi)\xi)E_1, 0, 0).
\end{array}$$
Since $\alpha(E_1, 0, 1, 0) \in (S')^{11}$, we have $(\tau\eta)\eta + \ov{y}y + (\tau\xi)\xi = 1$. Thus, we obtain $\alpha(E_1, 0,$ $ 1, 0) \times \alpha(E_1, 0, 1, 0) = 1/2{\mit\Phi}(0, E_1, 0, 0 )$, that is, $\alpha{\mit\Phi}(0, E_1, 0, 0)\alpha^{-1} = {\mit\Phi}(0, E_1, 0, 0)$. Since $\alpha \in Spin(12) \subset E_7$ satisfies $\alpha\tau\lambda = \tau\lambda\alpha$, we have also  $\alpha{\mit\Phi}(0, 0, E_1, 0)\alpha^{-1} = {\mit\Phi}(0, 0, E_1, 
\vspace{3mm}
0)$.

{\bf Proposition 4.4.} \qquad \qquad ${G_{12}}^{com} = Spin(12)$.
\vspace{2mm}
 
{\bf Proof}.  Now, let $\alpha \in {G_{12}}^{com}$. From $\alpha({\mit\Phi}_1, 0, 0, 0, 1, 0) = ({\mit\Phi}_1, 0, 0, 0, 1, 0)$ and $\alpha ({\mit\Phi}_1, 0, 0, 0, -1, 0) = ({\mit\Phi}_1, 0, 0, 0, -1, 0)$, we have $\alpha (0, 0, 0, 0, 1, 0)$  $ = (0, 0, 0, 0, 1, 0)$. Hence, since $\alpha \in {G_{12}}^{com} \subset E_8$, we see that $\alpha \in E_7$. We first show that $\kappa\alpha = \alpha\kappa$. Since ${G_{12}}^{com} \subset E_7$, it suffices to consider the actions on $\gP^C$. Since $\alpha \in {G_{12}}^{com}$ satisfies $\wti{\kappa}\alpha = \alpha\wti{\kappa}$, from 
$$
      \wti{\kappa}\alpha P =  \kappa\alpha P - \alpha P \quad \mbox{and}
       \quad 
      \alpha\wti{\kappa}P = \alpha\kappa P - \alpha P, \quad P \in \gP^C,$$
we have $\kappa\alpha = \alpha\kappa$. Next, we show that $\mu\alpha = \alpha\mu$.  
Again, from $\alpha ({\mit\Phi}_1, 0, 0, 0, 1, 0) = ({\mit\Phi}_1, 0, 0, 0, 1, 0)$ and $\alpha ({\mit\Phi}_1, 0, 0, 0, -1, 0) = ({\mit\Phi}_1, 0, 0, 0, 
 -1, 0)$, we have $\alpha ({\mit\Phi}_1, 0, 0, 0, 0, 0)$ $ = ({\mit\Phi}_1, 0, 0, 0, 0, 0)$. 
 Hence, since $\alpha \in E_7$, we have $\alpha{\mit\Phi}_1\alpha^{-1} = {\mit\Phi}_1$, that is, $\alpha{\mit\Phi}(0, E_1, 0, 0)\alpha^{-1}$ $ = {\mit\Phi}(0, E_1, 0, 0)$. Consequently
$$
\begin{array}{l}
   \alpha({\mit\Phi}(0, 0, E_1, 0), 0, 0, 0, 0, 1) = \alpha(-\wti{\mu}_\delta({\mit\Phi}(0, E_1, 0, 0), 0, 0, 0, 1, 0))
\vspace{1mm}\\
\quad
   = -\wti{\mu}_\delta\alpha({\it\Phi}(0, E_1, 0, 0), 0, 0, 0, 1, 0)
   = -\wti{\mu}_\delta({\it\Phi}(0, E_1, 0, 0), 0, 0, 0, 1, 0)
\vspace{1mm}\\
\quad
   = ({\mit\Phi}(0, 0, E_1, 0), 0, 0, 0, 0, 1).
\end{array}$$
Similarly, we have $\alpha({\mit\Phi}(0, 0, E_1, 0), 0, 0, 0, 0, -1) = ({\mit\Phi}(0, 0, E_1, 0), 0, 0, 0, 0, -1)$.\\
 Hence we have $\alpha({\mit\Phi}(0, 0, E_1, 0), 0, 0, 0, 0, 0) = ({\mit\Phi}(0, 0, E_1, 0), 0, 0, 0, 0, 0)$. Moreover, from $\alpha \in E_7$, we have $\alpha{\mit\Phi}(0, 0, E_1, 0)\alpha^{-1} = {\mit\Phi}(0, 0, E_1, 0)$. Hence put together with $\alpha{\mit\Phi}(0, E_1, 0, 0)\alpha^{-1} = {\mit\Phi}(0, E_1, 0, 0)$, we have $\alpha{\mit\Phi}(0, E_1, E_1, 0)\alpha^{-1} = {\mit\Phi}(0, E_1, E_1, 0)$, that is, $\alpha\mu\alpha^{-1} = \mu$. Thus,  we have $\mu\alpha = \alpha\mu$.    
Therefore, $\alpha \in (E_7)^{\kappa,\mu} = Spin(12)$.\\ 
\hspace*{3pt} Conversely, let $\alpha \in Spin (12)$. For $R \in {\gge_8}^C$, 
\begin{eqnarray*}
     \wti{\kappa}\alpha R \!\!\! &=& \!\!\! \mbox{[}(\kappa ,0, 0, -1, 0, 0 ), (\alpha{\mit\Phi}\alpha^{-1}, \alpha P ,\alpha Q, r, u, v)\mbox{]}
\vspace{1mm}\\
      \!\!\! &=& \!\!\! (\mbox{[}\kappa,\alpha{\mit\Phi}\alpha^{-1}\mbox{]}, \kappa\alpha P - \alpha P, \kappa\alpha Q + \alpha Q, 0, -2u, 2v)
\end{eqnarray*}
\noindent and
\begin{eqnarray*}
      \alpha\wti{\kappa}R \!\!\! &=& \!\!\! \alpha\mbox{[}((\kappa ,0, 0, -1, 0, 0), ({\mit\Phi}, P, Q, r, u, v)\mbox{]}
\vspace{1mm}\\
       \!\!\! &=& \!\!\! \mbox{[}\alpha(\kappa ,0, 0, -1, 0, 0), \alpha({\mit\Phi}, P, Q, r, u, v)\mbox{]}
\vspace{1mm}\\
   \!\!\! &=& \!\!\! (\mbox{[}\alpha\kappa\alpha^{-1},\alpha{\mit\Phi}\alpha^{-1}\mbox{]}, {\alpha\kappa\alpha^{-1}}(\alpha P) - 
     \alpha P, {\alpha\kappa\alpha^{-1}}(\alpha Q) + \alpha Q, 0, -2u, 2v).
\end{eqnarray*}
\noindent From $\kappa\alpha = \alpha\kappa$, we have $ \mbox{[}\alpha\kappa\alpha^{-1},\alpha{\mit\Phi}\alpha^{-1}\mbox{]} = \mbox{[}\kappa,\alpha{\mit\Phi}\alpha^{-1}\mbox{]}$. Thus, we have $\wti{\kappa}\alpha R = \alpha\wti{\kappa} R$, that is, $\wti{\kappa}\alpha = \alpha\wti{\kappa}$.
Next, from $\mu\alpha = \alpha\mu$ and Lemma 4.3, we have $\mu_1(\alpha{\mit\Phi}_1{\alpha^{-1}}){\mu_1}^{-1}\\
 = \alpha(\mu_1{\mit\Phi}_1{\mu_1}^{-1})\alpha^{-1} = \alpha{\mit\Phi}( 0, 0, E_1, 0 )\alpha^{-1} = {\mit\Phi}( 0, 0, E_1, 0 )$.
Hence, for $R = (\zeta{\mit\Phi}_1, P, 0,\\
 0, u, 0) \in (V^C)^{14}$, 
\begin{eqnarray*}
   \wti{\mu}_{\delta}\alpha R\!\!\!& = &\!\!\! \wti{\mu}_{\delta}(\zeta\alpha{\mit\Phi}_1\alpha^{-1}, \alpha P, 0, 0, u, 0)
\vspace{1mm}\\
    \!\!\! &=& \!\!\! ({\mit\Phi}(0, 0, -uE_1, 0), 0, i{\mu_1}\alpha P, 0, 0, -\zeta)
\end{eqnarray*}
\noindent and 
$$
\begin{array}{l}
     \alpha\wti{\mu}_{\delta}R = \alpha({\mit\Phi}(0, 0, -uE_1, 0), 0, i\mu_1 P, 0, 0, -\zeta)
\vspace{1mm}\\      
\qquad \quad \!\!
    = (\alpha{\mit\Phi}( 0, 0, -uE_1, 0)\alpha^{-1}, 0, i\alpha\mu_1 P, 0, 0, 
-\zeta)
\vspace{1mm}\\
\qquad \quad \!\!
   = ({\mit\Phi}(0, 0, -uE_1, 0), 0, i\alpha\mu_1 P, 0, 0, -\zeta).
\end{array}$$
\noindent Hence, from $\mu\alpha = \alpha\mu$, we have $\wti{\mu}_{\delta}\alpha R = \alpha\wti{\mu}_{\delta}R, R \in (V^C)^{14}$. From Lemma 4.3, we have $\alpha({\mit\Phi}_1, 0, 0, 0, 0, 0) = ({\mit\Phi}_1, 0, 0, 0, 0, 0)$. Moreover, since $\alpha \in E_7$, we have $\alpha(0, 0, 0, 0, 1, 0) = (0, 0, 0, 0, 1, 0)$ and $\alpha(0, 0, 0, 0, -1, 0) = (0, 0, 0, 0, -1, 0)$.
Hence, we have 
$ \alpha({\mit\Phi}_1, 0, 0, 0, 1, 0) = ({\mit\Phi}_1, 0, 0, 0, 1,0 )\,\,
\mbox{and}
    \,\,\alpha({\mit\Phi}_1, 0, 0, 0, -1, 0) = ({\mit\Phi}_1, 0, 0, 0,-1,0 ).$

Therefore, $\alpha \in {G_{12}}^{com}$. Thus, the proof of the proposition is completed.
\vspace{3mm}

{\bf Lemma 4.5}.\quad{\it The Lie algebras} ${\gg_{14}}^{com}$ and ${\gg_{13}}^{com} $ {\it of the groups} ${G_{14}}^{com}$ and \\
${G_{13}}^{com}$ {\it are given respectively by} 
$$ 
 \begin{array}{l}
 {\gg_{14}}^{com} = \{R \in \gg_{14}\,|\,\tau\wti{\lambda}(\ad R) = (\ad R)\tau\wti{\lambda} \}
\vspace{1mm}\\
\qquad  \quad \, = \Big\{\Bigl({\mit\Phi}\Bigl(D + \pmatrix{0 & 0 & 0 \cr
                                0 & 0 & d_1 \cr
                                0 & - \ov{d}_1 & 0}^{\!\!\sim} + 
                       i\pmatrix{\epsilon_1 & 0 & 0 \cr
                                0 & \epsilon_2 & t_1 \cr
                                0 & \ov{t}_1 & \epsilon_3}^{\!\!\sim},
                       \pmatrix{0 & 0 & 0 \cr
                                0 & \rho_2 & a_1 \cr
                                0 & \ov{a}_1 & \rho_3},
\end{array}$$
%\vspace{1mm}\\
%\qquad \qquad \quad \,
$$
\begin{array}{l}
\qquad 
        -\tau \pmatrix{0 & 0 & 0 \cr
                                0 & \rho_2 & a_1 \cr
                                0 & \ov{a}_1 & \rho_3}, \nu \Bigl), 
                 \Bigl(\pmatrix{0 & 0 & 0 \cr
                   0 & \zeta_2 & z_1 \cr
                   0 & \ov{z}_1 & \zeta_3},
          \pmatrix{\zeta_1 & 0 & 0 \cr
                   0 & 0 & 0 \cr
                   0 & 0 & 0},
                  0, \zeta \Bigl),
\vspace{1mm}\\
\qquad          
      -\tau\lambda\Bigl(\pmatrix{0 & 0 & 0 \cr
                   0 & \zeta_2 & z_1 \cr
                   0 & \ov{z}_1 & \zeta_3},
          \pmatrix{\zeta_1 & 0 & 0 \cr
                   0 & 0 & 0 \cr
                   0 & 0 & 0},
                  0, \zeta \Bigl), r, 0, 0 \Bigl)
\vspace{1mm}\\ 
\qquad   |\, D \in \so(8)
, \epsilon_i \in \R,\, \rho_i, \zeta_i, \zeta \in C,\, \nu, r \in i\R,\,  \epsilon_1 + \epsilon_2 + \epsilon_3 = 0, 
\vspace{1mm}\\
\qquad 
    i\epsilon_1 + \dfrac{2}{3}\nu + 2r = 0,\, d_1, t_1 \in \gC, \, a_1, z_1 \in \gC^C \Big\}. 
\end{array}$$    
$$ 
\begin{array}{l}
  {\gg_{13}}^{com} = \{R \in {\gg_{14}}^{com} \,| \, (\ad R)({\mit\Phi}_1, 0, 0, 0, 1, 0 ) = 0 \} 
\vspace{1mm}\\
\qquad \quad \,
    = \Big\{\Bigl({\mit\Phi}\Bigl(D + \pmatrix{0 & 0 & 0 \cr
                                0 & 0 & d_1 \cr
                                0 & - \ov{d}_1 & 0}^{\!\!\sim} + 
                       i\pmatrix{\epsilon_1 & 0 & 0 \cr
                                0 & \epsilon_2 & t_1 \cr
                                0 & \ov{t}_1 & \epsilon_3}^{\!\!\sim},
                       \pmatrix{0 & 0 & 0 \cr
                                0 & \rho_2 & a_1 \cr
                                0 & \ov{a}_1 & \rho_3},    
\vspace{1mm}\\
\qquad  \qquad  
            -\tau\pmatrix{0 & 0 & 0 \cr
                          0 & \rho_2 & a_1 \cr
                          0 & \ov{a}_1 & \rho_3}, \nu \bigl),

   \Bigl(\pmatrix{0 & 0 & 0 \cr
                   0 & \zeta_2 & z_1 \cr
                   0 & \ov{z}_1 & -\tau\zeta_2},
          \pmatrix{\zeta_1 & 0 & 0 \cr
                   0 & 0 & 0 \cr
                   0 & 0 & 0},
                  0, \tau\zeta_1 \Bigl),
\vspace{1mm}\\
\qquad          
      -\tau\lambda\Bigl(\pmatrix{0 & 0 & 0 \cr
                   0 & \zeta_2 & z_1 \cr
                   0 & \ov{z}_1 & -\tau\zeta_2},
          \pmatrix{\zeta_1 & 0 & 0 \cr
                   0 & 0 & 0 \cr
                   0 & 0 & 0},
                  0, \tau\zeta_1 \Bigl), 0, 0, 0 \Bigl)
\vspace{1mm}\\
\qquad \qquad 
  |\, D \in \so(8)
, \epsilon_i \in \R,\, \rho_i, \zeta_i \in C,\, \nu \in i\R,  \epsilon_1 + \epsilon_2 + \epsilon_3 = 0, \\
\qquad \qquad 
   i\epsilon_1 + \dfrac{2}{3}\nu = 0,\, d_1, t_1, z_1 \in \gC,\,  a_1 \in \gC^C \Big\}. 
\end{array}$$
\noindent {\it  In particular, we have} 
$$
      \dim({\gg_{14}}^{com}) = 28 + 63 = 91, \quad \dim({\gg_{13}}^{com}) = 28 + 50 = 78. $$

{\bf Lemma 4.6.}\,(1) {\it For $a \in \gC$,  we define a $C$-linear transformation $\epsilon_{13}(a)$ of ${\gge_8}^C$ by}
$$
     \epsilon_{13}(a) = \mbox{exp}(\ad(0, (F_1(a), 0, 0, 0), (0, F_1(a), 0, 0), 0, 0, 0)). $$
\noindent {\it Then}, $\epsilon_{13}(a) \in {G_{13}}^{com}$ ({\it Lemma} 4.5). {\it The action of $\epsilon_{13}(a)$ on $V^{13}$ is given by} 
$$
\begin{array}{l}
    \epsilon_{13}(a) ({\mit\Phi}(0, \zeta E_1, 0, 0), (\xi E_1, \eta E_2 - \tau\eta E_3 + F_1(y), \tau\xi, 0), 0, 0, -\zeta, 0)
\vspace{1.5mm}\\
\qquad
    = ({\mit\Phi}(0, \zeta'E_1, 0, 0), (\xi'E_1, {\eta}'E_2 - \tau{\eta}'E_3 + F_1(y'), \tau\xi', 0), 0, 0, -\zeta', 0),
\end{array}$$
$$
     \left \{ \begin{array}{l}
       \zeta' = \zeta\cos|a|
                 - \dfrac{(a, y)}{2|a|}\sin|a| 
\vspace{1mm}\\
       \xi' = \xi 
\vspace{1mm}\\
       \eta' = \eta 
\vspace{1mm}\\
       y' = y + \dfrac{2\zeta a}{|a|}\sin|a| - \dfrac{2(a, y)a}{|a|^2}\sin^2\dfrac{|a|}{2}. 
\end{array} \right. $$    

(2) {\it For $t \in \R$, we define a $C$-linear transformation $\theta_{13}(t)$ of ${\gge_8}^C$ by}
$$
     \theta_{13}(t) = \exp(\ad(0, (0, -tE_1, 0, - t), (tE_1, 0, t, 0), 0, 0, 0)). $$
\noindent {\it Then}, $\theta_{13}(t) \in {G_{13}}^{com}$ ({\it Lemma} 4.5). {\it The action of $\theta_{13}(t)$ on $V^{13}$ is given by} 
$$
\begin{array}{l}
   \theta_{13}(t)({\mit\Phi}(0, \zeta E_1, 0, 0), (\xi E_1, \eta E_2 - \tau\eta E_3 + F_1(y), \tau\xi, 0 ), 0, 0, -\zeta, 0)
\vspace{1.5mm}\\
\qquad
    = ({\mit\Phi}(0, {\zeta}'E_1, 0, 0), ({\xi}'E_1, {\eta}'E_2 - \tau{\eta}' E_3 + F_1(y'), \tau{\xi}', 0), 0, 0, -\zeta', 0),
\end{array}$$
$$
     \left \{ \begin{array}{l}
\vspace{1mm}
       \zeta' = \zeta\cos t - \dfrac{1}{4}(\tau\xi + \xi)\sin t 
\vspace{1mm}\\
       \xi' = \dfrac{1}{2}(\xi - \tau\xi)
                 + \dfrac{1}{2}(\xi + \tau\xi)\cos t
                 + 2\zeta\sin t 
\vspace{1mm}\\
       \eta' = \eta
\vspace{1mm}\\
       y' = y . 
 \end{array} \right. $$
\vspace{1mm}    

{\bf Lemma 4.7.} \qquad \qquad ${G_{13}}^{com}/{G_{12}}^{com} \simeq S^{12}$.
\vspace{1mm}

\noindent {\it In partiular,} ${G_{13}}^{com}$ {\it is connected}.
\vspace{2mm}

{\bf Proof.} Let $S^{12} = \{ R \in V^{13} \, | \, (R, R)_{\mu} = 1 \}$. The group ${G_{13}}^{com}$ acts on $(S^C)^{12}$. We shall show that this action is transitive. To prove this, it suffices to show that any $R \in S^{12}$ can be transformed to $1/2({\mit\Phi}_1, 0, 0, 0, -1, 0) \in S^{12}$. Now for a given 
$$
   R = ({\mit\Phi}(0, \zeta E_1, 0, 0), (\xi E_1, \eta E_2 - \tau\eta E_3 + F_1(y), \tau\xi, 0), 0, 0, -\zeta, 0) \in S^{12}, $$
\noindent choose $a \in \gC$ such that $|a| = \pi/2, (a, y) = 0$. Operate $\epsilon_{13}(a) \in {G_{13}}^{com}$ (Lemma 4.6.(1)) on $R$. Then, we have
$$
   \epsilon_{13}(a)R = (0, (\xi E_1, \eta E_2 - \tau\eta E_3 + F_1(y'), \tau\xi, 0 ),0 ,0 , 0, 0 ) = R_1 \in (S')^{11} \subset S^{12},$$
\noindent where $(S')^{11} = \{ R \in (V')^{12} \, | \, (R, R)_{\mu} = 1 \}$. Here, since the group $Spin(12)(\subset {G_{13}}^{com})$ acts transitively on $S^{11} = \{ P \in V^{12} \,| \, (P, P)_{\mu} = 1 \}$, there exists $\beta \in Spin(12)$ such that $\beta P = (0, E_1, 0, 1)$ for any $ P \in S^{11}$. Hence we have
\begin{eqnarray*}
                \beta R_1\!\!\! &=&\!\!\! \beta(0, P', 0, 0, 0, 0) = (0, \beta P', 0, 0, 0, 0)
  \vspace{1mm}\\
                           &=&\!\!\!(0, \beta\mu P, 0, 0, 0, 0) = (0, \mu\beta P, 0, 0, 0, 0)
  \vspace{1mm}\\
                           &=&\!\!\!(0, \mu(0, E_1, 0, 1), 0, 0, 0, 0) = (0, ( E_1, 0, 1, 0), 0, 0, 0, 0) = R_2 \in(S')^{11} ,
\end{eqnarray*}
\noindent where $ P \in S^{11}$. \\
\noindent Finally, operate $\theta_{13}(-\pi/2) \in {G_{13}}^{com}$ (Lemma 4.6.(2)) on $R_2$. Then, we have 
$$
        \theta_{13}(-\pi/2)R_2 = \dfrac{1}{2}({\mit\Phi}_1, 0, 0, 0, - 1, 0). $$
\noindent This shows the transitivity. The isotropy subgroup at $1/2({\mit\Phi}_1, 0, 0, 0, - 1, 0)$ of ${G_{13}}^{com}$ is obviously ${G_{12}}^{com}$. Thus, we have the homeomorphism ${G_{13}}^{com}/{G_{12}}^{com} $ $\simeq S^{12}$.
\vspace{3mm}

{\bf Proposition 4.8.} \qquad \qquad ${G_{13}}^{com} \cong Spin(13).$
\vspace{2mm}

{\bf Proof}. Since the group ${G_{13}}^{com}$ is connected (Lemma 4.7), we can define a homomorphism $\pi : {G_{13}}^{com} \to SO(13) = SO(V^{13})$ by 
$$
               \pi(\alpha) = \alpha|V^{13}.$$ 
\noindent $\Ker\,\pi = \{1, \sigma\} = \Z_2$. Since $\dim({\gg_{13}}^{com}) = 78 $ (Lemma 4.5) $= \dim(\so(13))$, $\pi$ is onto. Hence,  ${G_{13}}^{com}/\Z_2 \cong SO(13)$. Therefore, ${G_{13}}^{com}$ is isomorphic to $Spin(13)$ as a double covering group of $SO(13) = SO(V^{13})$.
\vspace{3mm}

{\bf Proposition 4.9.} \qquad \qquad ${G_{14}}^{com} \cong Spin(14)$.
\vspace{2mm}

{\bf Proof}. Since the group ${G_{14}}^{com}$ acts on $V^{14}$ and ${G_{14}}^{com}$ is connected(Proposition 4.2), we can define a homomorphism $ \pi : {G_{14}}^{com} \to SO(14) = SO(V^{14})$ by
$$   
              \pi(\alpha) = \alpha|V^{14}. $$
\noindent $\Ker\,\pi = \{1, \sigma\} = \Z_2$. Since $\dim({\gg_{14}}^{com}) = 91$ (Lemma 4.5) $= \dim(\so(14))$, $\pi$ is onto. Hence, ${G_{14}}^{com}/\Z_2 \cong SO(14)$. Therefore, ${G_{14}}^{com}$ is isomorphic to $Spin(14)$ as a double covering group of $SO(14) = SO(V^{14})$.
\vspace{2mm}

Now, we shall consider the following group
$$
\begin{array}{l}
   ((Spin(13))^{\sigma'})_{(0,F_1(y),0,0)^-}
\vspace{1mm}\\
\quad
   = \Big\{\alpha \in (Spin(13))^{\sigma'} \, \Big| \begin{array}{l} \alpha(0, (0, F_1(y), 0, 0), 0, 0, 0, 0) \\
= (0, (0, F_1(y), 0, 0), 0, 0, 0, 0) 
\end{array} \mbox{for all}\; y \in \gC
 \Big\}.
\end{array}$$
\vspace{-3mm}

{\bf Lemma 4.10.} {\it The Lie algebra $((\spin(13))^{\sigma'})_{(0,F_1(y),0,0)^-}$ of the group\\
 $((Spin(13))^{\sigma'})_{(0,F_1(y),0,0)^-}$ is given by}
$$
\begin{array}{l}
  ((\spin(13))^{\sigma'})_{(0,F_1(y),0,0)^-}\!= \! \{ R \in (\spin(13))^{\sigma'} | (\mbox{ad}R)(0, (0, F_1(y), 0, 0), 0, 0, 0, 0)\!= \! 0 \}
\end{array}$$
$$
\begin{array}{l}
           = \Big\{\Bigl({\mit\Phi}\Bigl( 
                      i\pmatrix{\epsilon_1 & 0 & 0 \cr
                                0 & \epsilon_2 & 0 \cr
                                0 & 0 & \epsilon_3}^{\!\!\sim},
                       \pmatrix{0 & 0 & 0 \cr
                                0 & \rho_2 & 0 \cr
                                0 & 0 & \rho_3},
                  -\tau\pmatrix{0 & 0 & 0 \cr
                                0 & \rho_2 & 0 \cr
                                0 & 0 & \rho_3}, \nu \Bigl), 
\end{array}$$
$$
\begin{array}{l}
\qquad
    \Bigl(\pmatrix{0 & 0 & 0 \cr
                   0 & \zeta_2 & 0 \cr
                   0 & 0 & -\tau\zeta_2},
          \pmatrix{\zeta_1 & 0 & 0 \cr
                   0 & 0 & 0 \cr
                   0 & 0 & 0},
                             0, \tau\zeta_1 \Bigl),
          -\tau\lambda\Bigl(\pmatrix{0 & 0 & 0 \cr                            
                                     0 & \zeta_2 & 0 \cr
                                     0 & 0 & -\tau\zeta_2},
\vspace{1mm}\\
\qquad 
          \pmatrix{\zeta_1 & 0 & 0 \cr
                   0 & 0 & 0 \cr
                   0 & 0 & 0},
                              0, \tau\zeta_1 \Bigl), 0, 0, 0 \Bigl)  
\vspace{1mm}\\
\qquad 
   |\, \epsilon_i \in \R,\, \rho_i, \zeta_i, \in C,\, \nu \in i\R,  \epsilon_1 + \epsilon_2 + \epsilon_3 = 0,\, i\epsilon_1 + \dfrac{2}{3}\nu = 0 \Big \},
\end{array}$$    
\noindent {\it In particular, we have} 
$$
    \dim(((\spin(13))^{\sigma'})_{(0,F_1(y),0,0)^-}) = 10. $$

{\bf Lemma 4.11.} \qquad \quad $((Spin(13))^{\sigma'})_{(0,F_1(y),0,0)^-}/Spin(4) \simeq S^4$.
\vspace{1mm}

\noindent {\it In particular}, $((Spin(13))^{\sigma'})_{(0,F_1(y),0,0)^-}$ {\it is connected.}
\vspace{2mm}

{\bf Proof}. We define a 5 dimensional $\R$-vector spaces $W^5$ by
$$
\begin{array}{l}
   W^5  = \{R \in V^{13} \, | \, \sigma'R = R \} 
\vspace{1mm}\\
\quad \;\;\;            
      = \{ R = ({\mit\Phi}(0, \zeta E_1, 0, 0), (\xi E_1, \eta E_2 - \tau\eta E_3, \tau\xi, 0), 0, 0, -\zeta, 0)
        \, |\, \zeta \in \R, \xi, \eta \in C \}
\end{array}$$
with the norm
$$ 
   (R, R)_{\mu} = \dfrac{1}{30}B_8(\wti{\mu}_{\delta}R, R) = 4\zeta^2 + (\tau\eta)\eta + (\tau\xi)\xi. $$
Then, $S^4 = \{ R \in W^5\,|$ $\, (R, R)_{\mu} = 1 \}$ is a 4 dimensional sphere. The group \\
$((Spin(13))^{\sigma'})_{(0,F_1(y),0,0)^-}$ acts on $S^4$. We shall show that this action is transitive. To prove this, it suffices to show that any $R \in S^4$ can be transformed to $1/2({\mit\Phi}_1, 0, 0, 0, -1, 0) \in S^4$ under the action of $((Spin(13))^{\sigma'})_{(0, F_1(y), 0, 0)^-}$. 
Now, for a given
$$
 R = ({\mit\Phi}(0, \zeta E_1, 0, 0), (\xi E_1, \eta E_2 - \tau\eta E_3 , \tau\xi, 0), 0, 0, -\zeta, 0) \in S^4, $$
choose $t \in \R, 0 \le t < \pi$ such that $\tan t = \dfrac{4\zeta}{\xi + \tau\xi}$ (if $\xi + \tau\xi$ = 0, let $t = \pi/2$). Operate $\theta_{13}(t) \in ((Spin(13))^{\sigma'})_{(0,F_1(y),0,0)^-}$ (Lemmas 4.6.(2), 4.10) on $R$. Then, we have
$$ 
     \theta_{13}(t)R = (0, ({\xi'} E_1, \eta E_2 - \tau\eta E_3, \tau\xi', 0), 0, 0, 0, 0) = R_1 \in S^3 \subset S^4.$$
Since the group $((Spin(12))^{\sigma'})_{(0,F_1(y),0,0)} (\subset ((Spin(13))^{\sigma'})_{(0,F_1(y),0,0)^-})$ acts transitively on $S^3$ (Lemma 3.14), there exists $\beta \in ((Spin(12))^{\sigma'})_{(0,F_1(y),0,0)}$ such that 
$$
    \beta R_1 = (0, (E_1, 0, 1, 0), 0, 0, 0, 0) = R_2 \in S^3. $$
Finally, operate $\theta_{13}(-\pi/2) \in ((Spin(13))^{\sigma'})_{(0,F_1(y),0,0)^-}$ on $R_2$. Then, we have 
$$ 
   \theta_{13}(-\pi/2)R_2 = \dfrac{1}{2}({\mit\Phi}_1, 0, 0, 0, -1, 0). $$
This shows the transitivity. The isotropy subgroup at $1/2({\mit\Phi}_1, 0, 0, 0, -1, 0)$ of \\
 $((Spin(13))^{\sigma'})_{(0,F_1(y),0,0)^-}$ is $((Spin(12))^{\sigma'})_{(0,F_1(y),0,0)}$ (Lemma 4.7) $ = Spin(4)$. Thus, we have the homeomorphism $((Spin(13))^{\sigma'})_{(0,F_1(y),0,0)^-}/Spin(4) \simeq S^4$.
\vspace{3mm}

{\bf Proposition 4.12}. \quad \qquad $((Spin(13))^{\sigma'})_{(0,F_1(y),0,0)^-} \cong Spin(5)$.
\vspace{2mm}

{\bf Proof}. Since $((Spin(13))^{\sigma'})_{(0,F_1(y),0,0)^-}$ is connected (Lemma 4.11), we can define a homomorphism $\pi : ((Spin(13))^{\sigma'})_{(0,F_1(y),0,0)^-} \to SO(5) = SO(W^5)$ by
$$
   \pi(\alpha) = \alpha|W^5. $$
$\Ker\,\pi = \{1, \sigma \} = \Z_2$. Since $\dim(((\spin(13))^{\sigma'})_{(0,F_1(y),0,0)^-}) = 10$ (Lemma 4.10) $= \dim(\so(5))$, $\pi$ is onto. Hence, $((Spin(13))^{\sigma'})_{(0,F_1(y),0,0)^-}/\Z_2 \cong SO(5)$. Therefore, $((Spin(13))^{\sigma'})_{(0,F_1(y),0,0)^-}$ is isomorphic to $Spin(5)$ as a double covering group of $SO(5)$.
\vspace{3mm}

{\bf Lemma 4.13.} {\it The Lie algebra $(\spin(13))^{\sigma'}$ of the group $(Spin(13))^{\sigma'}$ is given by}
\begin{center}
$\begin{array}{l}
  (\spin(13))^{\sigma'}
\vspace{1mm}\\
\quad
           = \Big\{\Bigl({\mit\Phi}\Bigl(D +  
                      i\pmatrix{\epsilon_1 & 0 & 0 \cr
                                0 & \epsilon_2 & 0 \cr
                                0 & 0 & \epsilon_3}^{\!\!\sim},
                       \pmatrix{0 & 0 & 0 \cr
                                0 & \rho_2 & 0 \cr
                                0 & 0 & \rho_3},
                  -\tau\pmatrix{0 & 0 & 0 \cr
                                0 & \rho_2 & 0 \cr
                                0 & 0 & \rho_3}, \nu \Bigl), 
%\vspace{1mm}\\
%\qquad
\end{array}$ \end{center}
$$
\begin{array}{l}
\quad 
     \Bigl(\pmatrix{0 & 0 & 0 \cr
                   0 & \zeta_2 & 0 \cr
                   0 & 0 & -\tau\zeta_2},
          \pmatrix{\zeta_1 & 0 & 0 \cr
                   0 & 0 & 0 \cr
                   0 & 0 & 0},
                             0, \tau\zeta_1 \Bigl),
          -\tau\lambda\Bigl(\pmatrix{0 & 0 & 0 \cr                            
                                     0 & \zeta_2 & 0 \cr
                                     0 & 0 & -\tau\zeta_2},
\vspace{1mm}\\
\qquad 
          \pmatrix{\zeta_1 & 0 & 0 \cr
                   0 & 0 & 0 \cr
                   0 & 0 & 0},
                              0, \tau\zeta_1 \Bigl), 0, 0, 0 \Bigl)    
\vspace{1mm}\\
\quad 
 |\, D \in \so(8), \epsilon_i \in \R,\, \rho_i, \zeta_i \in C, \nu \in i\R,\,  \epsilon_1 + \epsilon_2 + \epsilon_3 = 0,\, i\epsilon_1 + \dfrac{2}{3}\nu = 0 \Big \}.
\end{array}$$    
\noindent {\it In particular, we have} 
$$
    \dim((\spin(13))^{\sigma'}) = 
28 + 10 = 38. $$
\vspace{-5mm}

Now, we shall determine the group structure of $(Spin(13))^{\sigma'}$.
\vspace{3mm}

{\bf Theorem 4.14.}\,\, $(Spin(13))^{\sigma'} \cong (Spin(5) \times Spin(8))/\Z_2, \Z_2 = \{(1, 1), (-1,  \sigma) \}$.
\vspace{2mm}

{\bf Proof.} Let $Spin(13) = {G_{13}}^{com}, Spin(5) = ((Spin(13))^{\sigma'})_{(0,F_1(y),0,0)^-}$ and \\
$Spin(8) = ((F_4)_{E_1})^{\sigma'} \subset ((E_6)_{E_1})^{\sigma'} \subset ((E_7)^{\kappa,\mu})^{\sigma'} \subset ({G_{13}}^{com})^{\sigma'}$ (Theorem 1.2, Propositions 4.4, 4.8). Now, we define a map $\varphi : Spin(5) \times Spin(8) \to (Spin(13))^{\sigma'}$ by
$$
   \varphi(\alpha, \beta) = \alpha\beta. $$
Then, $\varphi$ is well-defined : $\varphi(\alpha, \beta) \in (Spin(13))^{\sigma'}$. Since $[R_D, R_5] = 0$ for $R_D = ({\mit\Phi}(D, 0, 0, 0), $ $0, 0, 0, 0, 0) \in \spin(8), R_5 \in \spin(5) = ((\spin(13))^{\sigma'})_{(0,F_1(y),0,0)^-}$ (Proposition 4.12), we have $\alpha\beta = \beta\alpha$. Hence, $\varphi$ is a homomorphism. $\Ker\,\varphi = \{(1, 1), (-1,  \sigma) \}$  $ = \Z_2$. Since $(Spin(13))^{\sigma'}$ is connected and $\dim(\spin(5) \oplus \spin(8))\! = \!10 \mbox{(Lemma 4.10)} + 28 = 38 =  \dim((\spin(13))^{\sigma'})$ (Lemma 4.13), $\varphi$ is onto. Thus, we have the isomorphism $(Spin(5) \times Spin(8))/\Z_2 \cong ((Spin(13))^{\sigma'}$. 
\vspace{2mm}

Now, we shall consider the following group
$$
\begin{array}{l}
   ((Spin(14))^{\sigma'})_{(0,F_1(y),0,0)^-}
\vspace{1mm}\\
\quad
   = \Big\{\alpha \in ((Spin(14))^{\sigma'} \, \Big| \begin{array}{l} \alpha(0, (0, F_1(y), 0, 0), 0, 0, 0, 0) \\
= (0, (0, F_1(y), 0, 0), 0, 0, 0, 0) 
\end{array} \mbox{for all}\; y \in \gC
 \Big\}.
\end{array}$$
\vspace{-3mm}

{\bf Lemma 4.15.} {\it The Lie algebra $((\spin(14))^{\sigma'})_{(0,F_1(y),0,0)^-}$ of the group\\
 $((Spin(14))^{\sigma'})_{(0,F_1(y),0,0)^-}$ is given by}
$$
\begin{array}{l}
  ((\spin(14))^{\sigma'})_{(0,F_1(y),0,0)^-}\! = \!\{ R \in (\spin(14))^{\sigma'} |(\mbox{ad}R)(0, (0, F_1(y), 0, 0), 0, 0, 0, 0)\! = \!0 \}
\end{array}$$
$$
\begin{array}{l}
           = \Big\{\Bigl({\mit\Phi}\Bigl( 
                      i\pmatrix{\epsilon_1 & 0 & 0 \cr
                                0 & \epsilon_2 & 0 \cr
                                0 & 0 & \epsilon_3}^{\!\!\sim},
                       \pmatrix{0 & 0 & 0 \cr
                                0 & \rho_2 & 0 \cr
                                0 & 0 & \rho_3},
                  -\tau\pmatrix{0 & 0 & 0 \cr
                                0 & \rho_2 & 0 \cr
                                0 & 0 & \rho_3}, \nu \Bigl), 
\end{array}$$
$$
\begin{array}{l}
\qquad
     \Bigl(\pmatrix{0 & 0 & 0 \cr
                   0 & \zeta_2 & 0 \cr
                   0 & 0 & \zeta_3},
          \pmatrix{\zeta_1 & 0 & 0 \cr
                   0 & 0 & 0 \cr
                   0 & 0 & 0},
                             0, \zeta \Bigl),
          -\tau\lambda\Bigl(\pmatrix{0 & 0 & 0 \cr                            
                                     0 & \zeta_2 & 0 \cr
                                     0 & 0 & \zeta_3},
\vspace{1mm}\\
\qquad 
          \pmatrix{\zeta_1 & 0 & 0 \cr
                   0 & 0 & 0 \cr
                   0 & 0 & 0},
                              0, \zeta \Bigl), r, 0, 0 \Bigl)    
\vspace{1mm}\\
\qquad 
   |\, \epsilon_i \in \R,\, \rho_i, \zeta_i, \zeta \in C,\, \nu, r \in i\R, \, \epsilon_1 + \epsilon_2 + \epsilon_3 = 0, \,i\epsilon_1 + \dfrac{2}{3}\nu + 2r = 0 \Big \}.
\end{array}$$    
\noindent {\it In particular, we have} 
$$
    \dim(((\spin(14))^{\sigma'})_{(0,F_1(y),0,0)^-}) = 15. $$

{\bf Lemma 4.16.} {\it For $t \in \R$,  we define a $C$-linear transformation $\theta_{14}(t)$ of ${\gge_8}^C$ by}
$$
     \theta_{14}(t) = \exp(\mbox{ad}(0, (0, itE_1, 0, it), (itE_1, 0, it. 0), 0, 0, 0)). $$
\noindent {\it Then}, $\theta_{14}(t) \in ((Spin(14))^{\sigma'})_{(0,F_1(y),0,0)^{-}}$ ({\it Lemma} 4.15). {\it The action of $\theta_{14}(t)$ on $V^{14}$} \\  {\it is given by} 
$$
\begin{array}{l}
    \theta_{14}(t)({\mit\Phi}(0, \zeta E_1, 0, 0), (\xi E_1, \eta E_2 - \tau\eta E_3 + F_1(y), \tau\xi, 0), 0, 0, -\tau\zeta, 0)
\vspace{1.5mm}\\
\qquad
    = ({\mit\Phi}(0, \zeta'E_1, 0, 0), (\xi'E_1, {\eta}'E_2 - \tau{\eta}' E_3 + F_1(y'), \tau{\xi}', 0), 0, 0, -\tau{\zeta}', 0),
\end{array}$$
$$
     \left \{ \begin{array}{l}
       \zeta' = \dfrac{1}{2}(\zeta + \tau\zeta) + \dfrac{1}{2}(\zeta - \tau\zeta)\cos t - \dfrac{i}{4}(\xi + \tau\xi)\sin t 
\vspace{1mm}\\
       {\xi}' = \dfrac{1}{2}(\xi - \tau\xi) + \dfrac{1}{2}(\xi + \tau\xi)\cos t - i(\zeta - \tau\zeta)\sin t  
\vspace{1mm}\\
       \eta' = \eta 
\vspace{1mm}\\
       y = y'.
\end{array} \right. $$    

{\bf Lemma 4.17.} \qquad \quad $((Spin(14))^{\sigma'})_{(0,F_1(y),0,0)^-}/Spin(5) \simeq S^5$.
\vspace{1mm}

{\it In particular}, $((Spin(14))^{\sigma'})_{(0,F_1(y),0,0)^-}$ {\it is connected.}
\vspace{2mm}

{\bf Proof}. We define a 6 dimensional $\R$-vector space $W^6$ by
\begin{eqnarray*}
   W^6  \!\!\! &=& \!\!\! \{R \in V^{14} \, | \, \sigma'R = R \}
\vspace{1mm}\\
      \!\!\! &=& \!\!\! \{R = ({\mit\Phi}(0, \zeta E_1, 0, 0), (\xi E_1, \eta E_2 - \tau\eta E_3, \tau\xi, 0), 0, 0, -\tau\zeta, 0) \, | \, \zeta, \xi, \eta \in C \}
\end{eqnarray*}
with the norm
$$       
   (R, R)_\mu = \dfrac{1}{30}B_8(\wti{\mu}_\delta R, R) = 4(\tau\zeta)\zeta + (\tau\eta)\eta + (\tau\xi)\xi. $$
Then,  $S^5 = \{ R \in W^6 \, | \, (R, R)_{\mu} = 1 \}$ is a 5 dimensional sphere. The group \\
$((Spin(14))^{\sigma'})_{(0,F_1(y),0,0)^-}$ acts on $S^5$. We shall show that this action is transitive. To prove this, it suffices to show that any $R \in S^5$ can be transformed to $1/2(i{\mit\Phi}_1, 0, 0, 0, i, 0) \in S^5$ under the action of $((Spin(14))^{\sigma'})_{(0,F_1(y),0,0)^-}$. 
Now, for a given
$$
 R = ({\mit\Phi}(0, \zeta E_1, 0, 0), (\xi E_1, \eta E_2 - \tau\eta E_3 , \tau\xi, 0), 0, 0, -\tau\zeta, 0) \in S^5, $$
choose $t \in \R, 0 \le t < \pi$ such that $\tan t = -\dfrac{2i(\zeta - \tau\zeta)}{\xi + \tau\xi}$ (if $\xi + \tau\xi$ = 0, let $t = \pi/2$). Operate $\theta_{14}(t) \in ((Spin(14))^{\sigma'})_{(0,F_1(y),0,0)^-}$ (Lemmas 4.15, 4.16) on $R$. Then, we have
$$ 
      \theta_{14}(t)R = ({\mit\Phi}(0, ({\zeta'} E_1, 0, 0), (\xi'E_1, \eta E_2 - \tau\eta E_3, \tau\xi', 0), 0, 0, -\zeta', 0) = R_1 \in S^4 \subset S^5.$$
Since the group $((Spin(13))^{\sigma'})_{(0,F_1(y),0,0)^-} (\subset ((Spin(14))^{\sigma'})_{(0,F_1(y),0,0)^-})$ acts transitively on $S^4$ (Lemma 4.11), there exists $\beta \in ((Spin(13))^{\sigma'})_{(0,F_1(y),0,0)^-}$ such that 
$$
    \beta R_1 = \dfrac{1}{2} ({\mit\Phi}_1, 0, 0, 0, -1, 0) = R_2 \in S^3. $$
Moreover, operate $\theta_{14}(\pi/2)$ and $\alpha(\pi/4)$ (Lemma 3.13) in order, 
$$
    \theta_{14}(\pi/2)R_2 = (0, (-iE_1, 0, i, 0), 0, 0, 0, 0) = R_3, $$
and
$$
    \alpha(\pi/4)R_3 = (0, (E_1, 0, 1, 0), 0, 0, 0, 0) = R_4. $$
Finally, operate $\theta_{14}(-\pi/2) \in ((Spin(14))^{\sigma'})_{(0,F_1(y),0,0)^-}$ on $R_4$. Then, we have 
$$ 
  \theta_{14}(-\pi/2)R_4 = \dfrac{1}{2}(i{\mit\Phi}_1, 0, 0, 0, i, 0). $$
This shows the transitivity. The isotropy subgroup at $1/2(i{\mit\Phi}_1, 0, 0, 0, i, 0)$ of \\
$((Spin(14))^{\sigma'})_{(0,F_1(y),0,0)^-}$ is $((Spin(13))^{\sigma'})_{(0,F_1(y),0,0)^-}$  (Proposition 4.8) $ =  Spin(5)$. Thus, we have the homeomorphism $((Spin(14))^{\sigma'})_{(0,F_1(y),0,0)^-}/Spin(5) \simeq S^5$.
\vspace{3mm}

{\bf Proposition 4.18}. \quad \qquad $((Spin(14))^{\sigma'})_{(0,F_1(y),0,0)^-} \cong Spin(6)$.
\vspace{2mm}

{\bf Proof}. Since $((Spin(14))^{\sigma'})_{(0,F_1(y),0,0)^-}$ is connected (Lemma 4.17), we can define a homomorphism $\pi : ((Spin(14))^{\sigma'})_{(0,F_1(y),0,0)^-} \to SO(6) = SO(W^6)$ by
$$
   \pi(\alpha) = \alpha|W^6. $$
$\Ker\,\pi = \{1, \sigma \} = \Z_2$. Since $\dim(((\spin(14))^{\sigma'})_{(0,F_1(y),0,0)^-}) = 15$ (Lemma 4.15) $= \dim(\so(6))$, $\pi$ is onto. Hence, $((Spin(14))^{\sigma'})_{(0,F_1(y),0,0)^-}/\Z_2 \cong SO(6)$. Therefore, $((Spin(14))^{\sigma'})_{(0,F_1(y),0,0)^-}$ is isomorphic to $Spin(5)$ as a double covering group of $SO(6)$.
\vspace{3mm}

{\bf Lemma 4.19.} {\it The Lie algebra $(\spin(14))^{\sigma'}$ of the group $((Spin(14))^{\sigma'}$ is given by}
$$
\begin{array}{l}
  (\spin(14))^{\sigma'}
\vspace{1mm}\\
\qquad
           = \Big\{\Bigl({\mit\Phi}\Bigl(D +  
                      i\pmatrix{\epsilon_1 & 0 & 0 \cr
                                0 & \epsilon_2 & 0 \cr
                                0 & 0 & \epsilon_3}^{\!\!\sim},
                       \pmatrix{0 & 0 & 0 \cr
                                0 & \rho_2 & 0 \cr
                                0 & 0 & \rho_3},
                  -\tau\pmatrix{0 & 0 & 0 \cr
                                0 & \rho_2 & 0 \cr
                                0 & 0 & \rho_3}, \nu \Bigl),
\end{array}$$ 
%\vspace{1mm}\\
%\qquad
$$
\begin{array}{l}
\qquad
    \Bigl(\pmatrix{0 & 0 & 0 \cr
                   0 & \zeta_2 & 0 \cr
                   0 & 0 & \zeta_3},
          \pmatrix{\zeta_1 & 0 & 0 \cr
                   0 & 0 & 0 \cr
                   0 & 0 & 0},
                             0, \zeta \Bigl),
          -\tau\lambda\Bigl(\pmatrix{0 & 0 & 0 \cr                            
                                     0 & \zeta_2 & 0 \cr
                                     0 & 0 & \zeta_3},
\vspace{1mm}\\
\qquad 
          \pmatrix{\zeta_1 & 0 & 0 \cr
                   0 & 0 & 0 \cr
                   0 & 0 & 0},
                              0, \zeta \Bigl), r, 0, 0 \Bigl) 
\vspace{1mm}\\
\qquad 
   |\, D \in \so(8), \epsilon_i \in \R,\,\rho_i, \zeta_i, \zeta \in C,\, \nu \in i\R,  \epsilon_1 + \epsilon_2 + \epsilon_3 = 0,\, i\epsilon_1 + \dfrac{2}{3}\nu + 2r = 0 \Big \}.
\end{array}$$    
\noindent {\it In particular, we have} 
$$
    \dim((\spin(14))^{\sigma'}) = 
28 + 15 = 43. $$
\vspace{-5mm}

Now, we shall determine the group structure of $(Spin(14))^{\sigma'}$.
\vspace{3mm}

{\bf Theorem 4.20.} \enskip $(Spin(14))^{\sigma'} \cong (Spin(6) \times Spin(8))/\Z_2, \Z_2 = \{(1, 1), (-1,  \sigma) \}$.
\vspace{2mm}

{\bf Proof.} Let $Spin(14) = {G_{14}}^{com}, Spin(6) = ((Spin(14))^{\sigma'})_{(0,F_1(y),0,0)^-}$ and \\
$Spin(8) = ((F_4)_{E_1})^{\sigma'} \subset ((E_6)_{E_1})^{\sigma'} \subset ((E_7)^{\kappa,\mu})^{\sigma'} \subset({G_{13}}^{com})^{\sigma'} \subset ({G_{14}}^{com})^{\sigma'}$ (Theorem 1.2, Propositions 4.8, 4.9). Now, we define a map $\varphi : Spin(6) \times Spin(8) \to (Spin(14))^{\sigma'}$ by
$$
   \varphi(\alpha, \beta) = \alpha\beta. $$
Then, $\varphi$ is well-defined : $\varphi(\alpha, \beta) \in (Spin(14))^{\sigma'}$. Since $[R_D, R_6] = 0$ for $R_D = ({\mit\Phi}(D, 0, 0, 0), $ $0, 0, 0, 0, 0) \in \spin(8), R_6 \in \spin(6) = ((\spin(14))^{\sigma'})_{(0,F_1(y),0,0)^-}$(Propo- \\ sition 4.18), we have $\alpha\beta = \beta\alpha$. Hence, $\varphi$ is a homomorphism. $\Ker\,\varphi = \{(1, 1), (-1,  \sigma) \}$ $ = \Z_2$. Since $(Spin(14))^{\sigma'}$ is connected and $\dim(\spin(6) \oplus \spin(8)) = 15 \mbox{(Lemma 4.15)} + 28 = 43 = \dim((\spin(14))^{\sigma'})$ (Lemma 4.19), $\varphi$ is onto. Thus, we have the isomorphism $(Spin(6) \times Spin(8))/\Z_2 \cong ((Spin(14))^{\sigma'}$. 
\vspace{2mm}

{\bf Acknowledgment.} The author records here his warmest gratitude to Professor Ichiro Yokota, who recommended to write this paper and has constantly encouraged him.

\bigskip
\begin{flushright}
\begin{tabular}{l}
Toshikazu Miyashita \\
Tohbu High School \\
Agata, Tohbu, 389-0517, Japan \\
E-mail: serano@janis.or.jp \\
\end{tabular}
\end{flushright}

\end{document}